# An extension of valid syllogisms to valid categorical arguments, and a reduction of the latter to only Barbara, Darapti, Darii and Disamis

To the memory of my sister Cristina Popa, 1948-2018.

DAN CONSTANTIN RADULESCU

*Abstract*
*One presents a simple Set Theory Model (STM) of the valid categorical arguments (VCAs) - a proper superset of the valid (categorical) syllogisms (VS). The main STM initiator was George Boole, who worked with a "universe of discourse", U, which contains the pairwise complementary sets, or categorical terms, S,S'(non-S),P,P'(non-P),M,M'(non-M), and is thus partitioned into eight subsets: SPM:= S∩P∩M, S'PM,...,S'P'M'. In STM all superfluous syllogistic figures are disregarded, and both the positive terms, S,P,M, and the negative terms, S',P',M', are allowed to appear in the pairs of categorical premises (PCPs) and their entailed logical consequences (LCs). This increases the number of distinct P (resp. S) premises from the six formulable via only positive terms, to eight, and the number of distinct PCPs from the 36 appearing in the Classical Categorical Syllogistic (CCS), to the 64 appearing in the STM. Out of the latter 64 PCPs, only 32 PCPs entail LCs and thus generate VCAs. These VCAs, (and the VS), split into four types. Each type contains eight VCAs. Any VCA can be re-written, via the term relabelings, p:= P↔P', s:=S↔S', m:=M↔M', and their compositions, ps, pm, sm, psm, as any other VCA of the same type. Thus the VCAs Barbara, Darapti, Darii, Disamis, can be chosen as type representatives for both VCAs and VS: via one of the above term relabelings, any VCA or VS can be re-written, without changing their PCP or LC contents, as either a Barbara, Darapti, Darii, or Disamis. Besides the VCAs and their LCs, one discusses simple/biliteral VCA sorites, empty set constraints (ESC), "VCA distribution conservation", and other Rules of Valid Categorical Arguments (RofVCA) which are "VCA generalized versions" of some of the Rules of Valid Syllogisms (RofVS). Both STM and CCS follow, for PCP classification purposes, the convention that the P term has to appear in the firstly listed premise. One compares the CCS, which defines the VS as PCPs formulable via only positive terms, whose entailed LCs are restricted to only the statements A(S,P), E(S,P), I(S,P), O(S,P), with the STM, whose VCAs also entail LCs of these other formats: A(P,S), E(S',P'), I(S',P'), O(P,S).*

## 1. Introduction

The universal set U contains the S, P, M sets, and also contains S', P', M' – the complements in U of S, P, M - thus U contains $2^3$=8 partitioning subsets – they will be called just subsets - no other set but a partitioning subset of U will be called a subset. Obviously, each of the S, P, M, S', P', M' sets is partitioned only into four subsets out of the eight subsets of U. Both CCS and STM accept LCs that only can be derived from both premises "put together" - immediate inferences from any of the premises are not acceptable LCs. According to Aristotle – see Striker 2009, p.20, "A syllogism is an argument in which, certain things being posited, something other than what was laid down results by necessity because these things are so." The key word here is "other", but its interpretations in CCS and STM differ. To obtain "new knowledge", in CCS, one eliminates the middle term M, in order to get a statement which connects only the S and P terms. In STM one finds subsets of U, on which both premises "synergistically" bestow attributes. For example, in CCS, Barbara's LC is A(S,P); in STM, Barbara's PCP entails two, more precise, LCs: S=SPM:=S∩P∩M, and P'=S'P'M'. To simplify the notation, more often than not, one denotes the term intersections by term adjacencies, the empty set by a 0 (zero) instead of the usual ∅ symbol, and one denotes the union of (usually disjoint) sets by a + (plus) instead of ∪ (union). Thus, SPM, (resp. S'P'M'), is the only subset of S, (resp. P') possibly remaining non-empty after Barbara's premises are applied. The middle terms are parts of their respective LCs, but they can be "dropped" to recover the CCS Barbara's LC: A(S,P)=A(P',S'). The next section spells out in detail STM's "one subset LC paradigm" - any PCP entails an LC if and only if it pinpoints to a unique subset of U, about which it asserts that is either definitely not empty, or, that it is the only subset possibly not empty of one of the sets S,P,M,S',P',M' - all the other three subsets of that set being emptied by two universal premises. Note also that the STM LC paradigm automatically satisfies Aristotle's requirement for a syllogism's LC: one definitely needs both premises to single out just one subset of U. The same paradigm suggests, or even requires, that the middle term should remain, contrary to Aristotle's belief, a part of the LC: a unique subset of U could not be singled out without it belonging to either M or M'. The above clearly indicates that PCPs containing two particular premises do not entail any LC – those 16 PCPs will be split below into two subtypes, (4a) and (4b). There are still PCPs of another type, which, according to either STM's "one subset LC paradigm", or to the Aristotelian definition of a syllogism's LC, do not entail any LCs. Any PCP of this latter type, is made of a universal premise emptying subsets of M, and a particular premise placing set elements in subsets of M' (or the other way around). They will be classified below, in both CCS and STM, as being of two subtypes also: (5a) and (5b). Note that CCS does not mention this latter PCP type as a separate class of PCPs which do not entail LCs. Any LC derived from such "no synergy" premises amounts to just repeating one or both premises – these PCPs are uninteresting, except for counting how many such PCPs do exist. In STM all terms are interpreted as sets, the PCPs, and LCs, contain both positive and negative terms, and the validity of a syllogism is decided via set relationships implied by both of its premises. Below one presents three methods, (which are, one hopes, even more enlightening then the CCS' usual 3-circle Venn diagram representation of the PCP terms), for finding the LC of any PCP. From them, it results that there exist six PCP types - the two types just described above - which do not entail any LC, and PCPs of other four types which entail LCs. The latter four PCP types generate in STM, (resp. CCS), four VCA, (resp. VS), types, with the VS set being a proper subset of the VCA set. Via term relabelings, any VCA can be re-written, without a change of its premises' content, nor of its LC content, as any other VCA of the same type. This is similar to the syllogisms' direct (ostensive) reduction via contrapositions and obversions, but without any transposition, (or order switch, or metathesis), of the premises, and without any S↔P relettering or S↔P role reassignment – those will be discussed below. (If the syllogistic figures are not used, then the simple conversions are unnecessary and are not even mentioned.) Supposing, (based on the Dictum de omni et nullo), the validity of the "perfect" Figure 1 syllogisms, in CCS, the syllogism reduction proves the validity of the syllogisms in other Figures. In STM, the validity of each VCA is made obvious by its LC pinpointing to a unique subset of U.

Note that one can disregard the superfluous syllogistic figures by using a matrix of all the distinct PCPs, whose rows, (resp. columns), are the distinct P-premises, (resp. S-premises). This will "shrink" the sets of PCPs, such as Ferio/Festino/Ferison/Fresison, which in fact differ in name only, to just one representative – call it Ferio's PCP: E(M,P)I(M,S). The other examples, Celarent/Cesare, Camestres/Camenes, Darii/Datisi, Disamis/Dimaris, Felapton/Fesapo are reduced to, respectively, Celarent, Camestres, Darii, Disamis and Felapton. Thus instead of the CCS "Figures' protocol" to build PCPs out of two premises containing three terms, M,P,S, which gives a total of 64 PCP choices: four choices E,A,I,O for the 1$^{st}$ (P) premise, times four choices for the 2$^{nd}$ (S) premise, times another four choices for the figure arrangements of the terms in each of the 16 (four by four) PCPs, one uses an STM "matrix protocol" to build only the really distinct PCPs. Part of the "matrix protocol" is that given three arbitrary concrete terms, one of them is chosen as the middle term – appearing in both premises, (and, as part of the "formalization sub-protocol" is denoted by M), a $2^{nd}$ term is chosen to always be used in the firstly listed premise, (and is denoted by P), and the $3^{rd}$ term is chosen to always be used in the secondly listed premise, (and is denoted by S). One finds out that there are only six distinct P, (resp. S), premises formulable via using only positive terms (as CCS requires). The "matrix protocol" gets rid of uselessly permuting arguments in E(M,P)=E(P,M) and I(M,S)=I(S,M), to account for figures, and one acknowledges only E(M,P) and I(M,S) as distinct P- premises, instead of considering that E(M,P), E(P,M), I(M,S), I(S,M) are all distinct. The fact that A:=A(M,*) differs from A(*,M)=E(M',*)=:E', (resp. O:=O(M,*) differs from O(*,M)=I(M',*) =:I'), where *∈{S, P}, is recognized via this new notation. Using syllogistic figures to account for the difference between the E' and I' statements and the A and O statements is inefficient – but set theory was not around during Aristotle's time. One may also introduce the A':=A(M',*) and O':=O(M',*) statements, where *∈{S, P}. This way one gets to eight P, (resp. S), premises, forming a cube of opposition, split into two squares of opposition, E,A,I,O, and E',A',I',O'. The STM uses an eight by eight PCP matrix, where any PCP is obtained by choosing a premise from the P-cube, (and listing it firstly),  and a second premise from the S-cube; the positive and negative terms appear on an equal footing. Similarly, with the universal premises listed firstly as, E,A,E',I,O,I', one sees that, in the six by six matrix of all the PCPs formulable only via positive terms, the PCPs of the Figure 1 syllogisms, Barbara, AE', Celarent, EE', Darii, AI, Ferio, EI, will appear above the matrix' diagonal. Again, both CCS and STM use the convention that the firstly listed premise contains the P term: according to this convention EI means E(M,P)I(M,S), and IE means I(M,P)E(M,S) – a sub-diagonal PCP matrix element. The CCS' supplementary condition that P be the predicate of the conclusion is really not a restriction on the PCPs, but on how a VS is defined; nevertheless, the fact that CCS is interested only in VS leads to discarding any VCA, (and resp. any PCP), which is a not a VS, (resp. does not generate a VS). Ferio's LC is O(S,P), and IE's LC is O(P,S); both PCPs generate VCAs – but in CCS the IE:O(P,S) VCA will be discarded. One can see that the sub-diagonal elements are obtained from the above-diagonal elements via a "P↔S role reassignment" followed by a switch in the premises' order. This means that given three terms, once the middle term is agreed upon - and denoted by M, in fact it does not matter which one of the other two is chosen as P, (and used in the firstly listed premise), since, in the PCP matrix, the roles of P and S are switched in the above-diagonal matrix elements PCPs, compared to the below-diagonal PCPs. On the other side, (see Quine, p.106), via a metathesis, followed by a P↔S relettering, (unnecessary if one deals with concrete terms, for which the relettering is implicitly performed, as soon as the metathesis is performed, due to the convention that the firstly listed premise contains the term denoted by the letter P), one may write, (without changing its content nor its LC's content), a PCP such as I(M,P)E(M,S) whose LC is O(P,S), as E(M,P)I(M,S):O(S,P) – where the LC is written after the column sign. But at most one of the initial PCPs, EI and IE, could be sound; thus, the sound one, if any, may be written as a Ferio. This may be taken as the reason behind the CCS restricting, by definition, the VS to only having LCs of the A(S,P), O(S,P), E(S,P), I(S,P) formats: the VCAs having A(P,S) and O(P,S) as LCs may be rewritten, without any content change, as a VS. (The other example is Bocardo: there is no Bacordo since it is impossible that both OA and AO PCPs are sound, even if each one of them may be written in both formats. But CCS has both a Barbari, and a Bramantip, even if both having sound premises implies S=M=P. Also, Celarent and Camestres having (all four) sound premises implies S=P=0. If the PCPs, of these two existential import (ei) VCAs, EA:O(S,P) (Felapton), and AE:O(P,S) (Falepton?), have only sound premises, it results that M=0, which contradicts their possible LCs, O(S,P), resp. O(P,S) which are based on the ei condition M≠∅. Section 13 examines "how many sound VS or VCAs one may hope to construct out of three given terms, without imposing restrictions on the structure of the universal set U".)

To recap, in either the six by six PCP matrix used by CCS, or in the eight by eight PCP matrix used by STM, the S↔P role reassignment and a metathesis maps the above the diagonal part of the matrix onto the under the diagonal part of it - and vice versa - while leaving the diagonal PCPs of the matrix unchanged. A metathesis and the S↔P relettering, re-writes the under the diagonal PCP matrix elements as above the diagonal matrix elements, (and vice versa), without changing their content. From this point of view one may say that via a relabeling from the group G={e, m, p, s, mp, ms, ps, mps}, plus a metathesis and a P↔S relettering – the latter two necessary only if the initial syllogism's firstly listed premise was particular – any syllogism can be written as either a Barbara, a Darapti, or a Darii. But this way, given three terms – one of which was already designated as the middle term, one may end up with two Darii, instead of one Darii and one Disamis. Thus the question: what is preferable, two Darii, two Disamis, or one Darii and one Disamis? I prefer one Darii and one Disamis, and thus I maintain that there are four PCP (and VCA) types not three. Examples: All men are perfect, Some men are serious: Therefore Some serious (men - keep the middle term!) are perfect. Or, All men are serious, Some men are perfect: Therefore Some perfect (men) are serious. [Which is the same as Some serious (men) are perfect.] Note also that, via just one obversion, E(M,P)=A(M,P'), Celarent, (and resp. Ferio), may be written as a Barbara, (and resp. Darii), in the variables S,P',M. The action of each element of the group G can be expressed as a combination of obversions and contrapositions since there are no other set preserving transformations. (Except for the simple conversions, which, in STM, are simply not worth mentioning – since the syllogistic figures were discarded.) The term relabelings are determined by which term is designated as being either S or S', which one is labeled as either P or P', and which one is labeled either M or M', and none of these three choices and their compositions is influenced by logic. Thus, one may say that, up to term relabelings, all the VCAs, (resp. VS), of the same type are equivalent. (This applies to the subtypes of PCPs which do not entail any LCs, too: all the PCPs in each of the sets (4a), (4b), (5a), (5b), are, up to term relabelings, equivalent.) Any LC of a VCA or VS pinpoints to just one and only one of the eight subsets of U – thus the middle term is a rightful part of the LC, since either M or M' is part of the identification of that unique subset of U to which the LC is pointing to. Eliminating, i.e., most often just dropping, the middle term from the LC – a requirement originating from Aristotle – only weakens the LC – the weakened LC will then "fuzzily pinpoint" to two subsets of U – one partially identified by M, the other by M'. In STM, any LC consists of one of the assertions: (i) one subset remains possibly non empty but the other three subsets are empty, out of the four subsets partitioning either one of the "end (or outer) terms"/sets S,P,S',P', - this type of LC is entailed by the Barbara type (1) PCPs – made of two universal premises which empty four subsets – two from M, and two from M'; (ii) one subset remains possibly non empty but the other three subsets are empty, out of the four subsets partitioning either M or M' - this type of LC is entailed by the Darapti type (2) PCPs – made of two universal premises, emptying three subsets of either M or M'; (iii) one subset out of eight is definitely not empty - this type of LC is entailed by the Darii type (3a) PCPs, or by the Disamis type (3b) PCPs – made of one universal premise, plus one particular premise, emptying subsets from, and respectively placing set elements into, either M or M' - but not both. (Listing the P-premises as rows, and the S-premises as columns, and listing the universal premises firstly, each PCP of type (3a) will appear above the diagonal of the PCP matrix, and to each PCP of type (3a), corresponds a PCP of type (3b), simmetrically situated under the diagonal. The difference between the (3a) and (3b) VCAs is that the P term appears in the (3a) universal premises, while the S term appears in the (3b) universal premises.) For example, both AI (Darii's PCP), and IA (Disamis' PCP) pinpoint to the same precise LC: the SPM subset is not empty. From one PCP to the other one, the S and P terms incur "role reassignment" but the two PCPs, and thus VCAs, are compatible – they can be simultaneosly sound. But EI (Ferio's PCP), and IE (Fireo's PCP?!) contain contradictory premises. If, e.g., IE: O(P,S), is the sound VCA, (and, consequently, Ferio has false premises), then, via a metathesis and an S↔P relettering, Fireo (?!) becomes Ferio, and, in the new lettering, IE, (aka the old EI), may be forgotten since its premises are both false. Since the PCP matrix contains all possible pairings of premises, it results that given any three terms, once the middle term denotted by M is chosen, it does not matter which one of the other two terms is denoted as P and which one is chosen as S. An initial S↔P name switch just flips the PCP matrix about its diagonal. For each of the four types of PCPs entailing at least one LC, the pointing to a unique subset of U is key: if a PCP singles out one of the eight subsets of U, then one has an LC – if not, there is no LC. As noted, pinpointing to a unique subset of U also means specifying if that subset is part of M or M' – information which is promptly lost when the middle term is eliminated from the LC. The same applies to the simple VCA sorites with "biliteral premises", (see Carroll (1977), p. 280), which one considers in Section 6: an LC pinpoints to a unique subset from the partitioning of a larger universe of discourse made from all the terms appearing in the sorite. (Note that, e.g., Venn (1981), Carroll (1977), and even Boole (1854), also considered "multiliteral" sorites - having premises connecting more than two terms at a time – for example, All x is y or z, written in STM as x=x∩(y∪z)=:x(y+z) =xy+xz.)

Most logicians, past and presents, seem to prefer the middle term(s) elimination approach to syllogistics, initiated by Aristotle, (and pursued by Boole), while Jevons (1869), seems to have been the first and strongest proponent of an approach aimed at finding chains of inclusions whose final result, after substitutions, is an intersection which defines - as being the LC - just one partitioning subset of the universal set U. Using Barbara's premises, one may easily illustrate the two approaches – the elimination approach is based on an "intersection principle", (if one finds a set whose intersections with two complementary sets are both empty sets, then the initial set is itself empty, too); the decomposition/substitution approach is based, as Jevons puts it, on the "substitution of similars" and is "the true principle of reasoning". Consider Barbara: All M is P, All S is M - statements which translate to MP':=M∩P'=0, SM'=0.

(i) One infers that SP'=0, since SP'=SP'M+ SP'M'=0+0=0, from which one infers that All S is P. Thus, one has to find an intersection of terms, (SP'), whose intersections with both a set, (M), and its complement, (M'), are empty – then the initial intersection of terms is itself empty. George Boole's "elimination algebra" is not needed.

(ii) Barbara's premises may also be written as M=MP+MP'=MP, S=SM+SM'=SM. Instead of trying, à la Boole, to eliminate M from the above two "equations", one uses substitutions, à la Jevons: S=SM=SMP. One could also have written the premises as: P'=P'M+P'M'=P'M', M'=M'S'. Using substitutions one gets, P'= P'M'=P'M'S'. Thus one obtains two LCs: each of the sets S and P' is reduced to just one subset of U.

As a more complicated example, consider a 10-term/8 premise multiliteral sorite from Carroll (1986) pp. 287-292: d'n'm'=0, ka'c'=0, lem=0, dhk'=0, h'la'=0, hm'b'=0, a'bn=0, am'e=0. The universe of discourse U is generated by 10 terms, and it contains $2^{10}$=1024 subsets.

Using the elimination approach, Carroll obtains the LC c'el=0, which contains only terms which appear in the premises without their complements appearing, too – these terms are the "retinends" - they can not be eliminated; the rest of the terms are "eliminands" - they are not part of the sought LC. The LC may also be written as el=cel, which shows that the intersection el is reduced to, or included into, the $2^{10}/2^3=2^7$=128 subsets of U – into which the cel set may be decomposed. But using the substitution approach one obtains that el=elm'a'hbn'dkc which shows that el is in fact reduced to a unique subset of U. Thus, whenever it works, one prefers Jevons' substitution approach which provides more information; the middle term(s) will be dropped from the precise LCs only to show that the traditional LCs can be recovered. [Numbering the above premises/equations from 1 to 8, one can eliminate b using the equations 7 and 6, since both b and b' appear in those equations – and nowhere else: it results that a'hm'n=0. But, for example, the elimination of m from equation 3, is "barred" - as Carroll's puts it – by the fact that m' appears in three other equations,1, 6, and 8, instead of only one – therefore, the elimination of m has to be postponed until the three equations in m' are reduced to only one equation. One can eliminate d using the equations 1 and 4, and it results that hk'm'n'=0. From this and the previously obtained a'hm'n=0, one can eliminate n to obtain a'hk'm'=0. From this and equation 2 one eliminates k to obtain a'c'hm'=0. From the latter and equation 5 one eliminates h and obtains a'c'lm'=0. From the latter and equation 8 one eliminates "a" to obtain c'elm'=0. By now the equations 1,6, and 8 were all used, and one can eliminate m from equation 3 and the just obtained c'elm'=0, to finally obtain Carroll's LC: c'el=0, which contains only the retinends. Using Jevons' "decomposition/substitution principle", (instead of the "intersection principle"), one can write el=elm'=elma'=...= elm'a'hbn'dkc, where the equations 3,8,5,6,7,1,4,2 were used in this order to discard the null sets from the successive decompositions of el with respect to (m,m'), then (a, a'), then (h, h'), then (b,b'), etc.]

In Section 11, one defines the Domain of Applicability (DofA) of the RofVS as being the set of PCPs for which the CCS and the RofVS predict the same LCs; it results that on their DofA, the RofVS can replace the CCS. The RofVS' DofA turns out to be made of four PCPs whose LCs are either A(P,S) or O(P,S), and also of the eight PCPs which generate the eight Boolean VS, (aka the 15 "standard" VS when syllogistic figures are taken into account) - see Hurley (pp.290-291), Copi (pp.245-248). (See also Sections 6 and 11.) By comparison, the RofVCA can predict both existential import (ei) and non existential import (non-ei) LCs; their DofA contains all the 32 PCPs which entail LCs and ei LCs – see Section 12. (Darapti, Barbari, Bramantip, are examples of ei VCAs and ei VS; they include an ei condition on M, S, and P, respectively.) Note that one can simply ignore the RofVS, as Quine, p. 107, does.

In STM one discards the two RofVS which refer to the PCPs only, "the middle term has to be distributed in at least one premise", and, "two negative premises are not allowed". In fact, the above two RofVS, unnecessarily, but interestingly, postulate the existence of another two types of PCPs which do not entail any LC. Then the fact that the PCPs of types (4a), (4b), (5a), (5b) do not entail any LCs can be deduced as theorems. For example, Stebbing, pp. 86-92, uses the above two RofVS to prove that PCPs made of two particular premises do not entail any LCs. To be precise, PCPs formulable only via positive terms, and where both premises are either negative, or, where M is undistributed in both premises, do not entail any of the "four standard LCs", which, by definition, a VS should have: A(S,P), E(S,P), I(S,P), O(S,P). (See Hurley pp. 280-282.) But there are a few PCPs, E'E':=E(M',P)E(M',S) =A(P,M) A(S,M), EE=E(M,P) E(M,S), OE, EO, (for detailed explanations about the notation see next section), all formulable via only positive terms, which, in fact, generate perfectly valid VCAs: in the first PCP, M is nowhere distributed, the other three PCPs contain only negative premises; for each one of these four PCPs, the LC or ei LC is I(S',P'). Copying Carroll's "example strategy" of making the negative terms as intuitive as the positive ones, consider this PCP: All the desire of fighting is in my boys. All the desire of playing sports is in my boys. The middle term, my boys, is undistributed, but if the universe of discourse is my children and their desires, then my boys' complementary set is my girls. So, through obversions, the above PCP becomes one made of two negative premises: No desire of fighting is in my girls. No desire of playing sports is in my girls. Which, as well as the initial PCP, has this LC: My girls have no desire of fighting, nor of playing sports. If the set of my girls is non-empty, the ei LC is that there are Some (girls - mine), who have no desire or fighting, nor of playing sports. The above are examples of the E'E' and EE premises, having

I(S',P') as LC. As another example of the limited applicability of RofVS #1 and #2, consider a VCA similar to Darii, where the middle term is M' not M – one may call it Darii' - whose premises are A(M',P)I(M',S). Because M' appears instead of M, both premises are negative, (see Section 10, Ladd-Franklin, or DeMorgan), and the premises can not be formulated using only positive terms, but the LC exists – it is M'SP≠Ø – which becomes I(S,P) if M', (which is included in P), is dropped from the LC. If one is interested only in recovering the LCs with the middle term removed from the LCs, (the classical style syllogistic!), one can generalize the rest of the RofVS, #3 to #6, (which refer to the PCPs and their LCs together), to four Rules of VCA (RofVCA) such that these RofVCA will predict the LCs for all the VCAs and ei VCAs; any such LC will be one of the eight A,E,I,O statements in which no middle term appears. They can all be written as $E(S^*,P^*)$ and $I(S^*,P^*)$, where $P^* \in \{P,P'\}$, $S^* \in \{S,S'\}$.

As mentioned, there are two PCP types, each split into two PCP subtypes, which do not entail any LC: (4a) $I(M^*,P^*) I(M^*,S^*)$ and (4b) $I(M^*,P^*)I(M^{*'},S^*)$ - each subtype contains eight PCPs for a total of 16 PCPs made of two particular premises, and, (5a) $E(M^*,P^*)I(M^{*'},S^*)$ and (5b) $I(M^*,P^*)E(M^{*'},S^*)$ - each subtype contains eight PCPs for a total of 16 PCPs made of one universal and one particular premises, one "acting on" M and the other on M', where $M^* \in \{M, M'\}$, $P^* \in \{P, P'\}$, $S^* \in \{S, S'\}$ – a notation which will be used everywhere within this paper. (Note that $M^{*'}=M$ if $M^*=M'$, etc.)

One may note Lewis Carroll's opinion about CCS at the end of the 19th century, Carroll (1977), p.250: "...the ordinary textbooks of Formal Logic have elaborately discussed no less than nineteen different forms of Syllogisms – each with its own special and exasperating Rules, while the whole constitutes an almost useless machine, for practical purposes, many of the Conclusions being incomplete, and many quite legitimate forms being ignored", "As to syllogisms, I find that their nineteen forms, with about a score of others which the textbooks have ignored, can all be arranged under three forms, each with a very simple Rule of its own". (Carroll considers the types (3a), Darii, and (3b), Disamis, as one and the same type. The "19 forms of syllogisms" referred by Carroll were the "15 standard VS", plus the ei VS Bramantip, Darapti, Felapton and Fesapo.) Carroll (1977), p.240 or Carroll, (1958), p.173, gives three concrete examples of VCAs containing only negative premises, and asserts that "The theory that two negative premises prove nothing" is "another craze of 'The Logicians', fully as morbid as their dread of a negative Attribute" (i.e., the negative terms S',P',M'). In a similar way, Quine, 1982, p.107, dismisses many of the notions related to syllogisms, except the (3-circle Venn) diagrammatic check of their validity. Quine, in the short chapter on syllogisms, pp.102-108, starts with VS examples in only six moods, (each mood in a specific figure), AAA, EAE, AII, EIO, AOO, OAO, and uses "figures' equivalencies", a metathesis and the relettering S↔P to recast them into all the 15 standard non-ei VS of the CCS. This paper shows – see Section 7 - that all the 15 standard (non-ei) VS are equivalent to either AAA, (Barbara), AII, (Darii), or, IAI, (Disamis), and, via a metathesis, (i.e. a premises' order switch), and a relettering S↔P, they can be all written as either AAA, Barbara, or, AII, Darii. Another method of proving these equivalencies would have been to use, in reverse order, the above cited Quine's arguments, to "reduce" the 15 standard VS to the six moods used by Quine, and then add to those six moods a proof that AAA and EAE are equivalent, and that AII, EIO, AOO, OAO may all be written as Darii or Disamis. For example, EAE, Celarent/Cesare, means E(M,P)A(S,M):E(S,P)= A(M,P')A(S,M): A(S,P'), where one writes the LC after the column sign. This shows that, by replacing P with (P')', Celarent/Cesare may be written as Barbara. For a more concrete illustration, the following Celarent and Barbara also say exactly the same thing: No M are boys, All S are M:No S are boys, and, All M are girls, All S are M: All S are girls. Thus, whatever a VCA of type (1) expresses, can also be expressed via an equivalent Barbara. In general, any non-ei VCA built using three given terms, out of which one has been chosen to serve as the middle term, is equivalent, via a relabeling, to either Barbara, Darapti, Darii, or Disamis, but, via a metathesis, any Darii may be written as a Disamis and vice versa. In STM, the P-term might end-up being the subject of the LC, i.e., LCs of the type A(P,S) and O(P,S) are VCA acceptable, even if they are not VS acceptable. The question of how many sound VCAs one may be able to build out of three given terms is discussed in Section 13.

Throughout, one stresses the STM's "one subset LC paradigm", and one insists on keeping the middle term as part of the precise LC.

The rest of the sections bring the necessary details about exact LCs, about LCs from which the middle term was eliminated, about biliteral sorites, about how many distinct VCAs one may construct from three given terms without imposing restrictions on the structure of the universal set U, about the RofVS and the RofVCA "theories", about the fact that all the VCAs (and VS) of the same type are equivalent, etc. (The equivalencies of the non-ei VCAs and those of the ei VCAs, are precisely stated at the end of Section 7. For example, if one refers only to the non-ei VCAs, whose LCs are the precise ones – i.e., LCs out of which the middle term was not eliminated, then one may say, in short, that any such VCA is equivalent to either Barbara (whose two LCs are S=SPM, and P'=S'P'M') , the non-ei Darapti, Darii, or Datisi. If the ei VCAs are included, stating the equivalencies is a little bit more complicated: see the end of Section 7.)

## 2. Summary of STM – an algebraic method for finding the LCs

The STM discards, but CCS uses, syllogistic figures. Moreover, in CCS only the statements A(S,P), E(S,P), I(S,P), O(S,P) are accepted as LCs, and only positive terms are accepted in its premises (and its LCs – as the above four LC accepted formats show). In STM, both positive and negative terms are accepted in PCPs and LCs. The positive only terms constraint implies that, contrary to what the syllogistic figures suggest, in CCS there are only six distinct P-premises and only six distinct S-premises, and thus there are only 36 distinct PCPs formulable via only positive terms - instead of the counted 64 PCPs - when syllogistic figures are used. As a result, for this one syllogism, "No M is P, Some M is S; therefore Some S is not P", the CCS has four names, Ferio/Festino/Ferison/Fresison, but one may use just the name Ferio, irrespective of the syllogistic figure; Darii/Datisi will be most often called Darii, etc. The four P-premises, A(M,P):=All M is P=E(M,P'), E(M,P):=No M is P, I(M,P):=Some M is P, O(M,P):= Some M is not P=I(M,P'), form a square of opposition containing M and P and may be denoted by A,E,I,O. The other two P-premises formulable using only the positive terms P and M, are A(P,M):=All P is M, and O(P,M):= Some P is not M. One notices, that in fact, A(P,M), (resp. O(P,M)), "act upon" M', since they empty, (resp. lay set elements in), the intersection PM':=P∩M'. Moreover, one can see that they act on M' exactly how E(M,P) and I(M,P) were acting on M: A(P,M)= E(M',P)=:E' and O(P,M)=I(M',P)=:I'. As a result, one may form another square of opposition, denoted A',E',I',O', which contains M' and P, by adding two more P-premises: A(M',P)= E(M',P')=:A' and O(M',P)=I(M',P')=:O'. The eight P-premises form now a cube of opposition, whose universal, (resp. particular), statements may be also written as E(M*,P*), (resp. I(M*,P*)), where M*∈{M, M'}, P*∈{P, P'}. Similarly, the eight S-premises will form another cube of opposition, E(M*,S*), I(M*,S*), where M* ∈{M, M'}, S* ∈{S, S'}. By choosing a premise from each cube, in STM, one gets 64 distinct PCPs, where the positive and negative terms appear "on an equal footing". Thus, one comes up with two compact notations. In one notation, any of the four PCP types which entail LCs, and any of the four PCP subtypes which do not entail any LC, may be represented by a single formula, e.g., the eight type (1) – Barbara – PCPs, may be written as E(M*,P*)E(M*',S*), where M*∈{M, M'}, P*∈{P, P'}, S*∈{S, S'}. And, by following the convention to firstly list the P-premise, one can use a 2-letter shorthand notation to write each PCP without explicitly listing its terms: for example, the 2-letter PCP A'A, would mean A(M',P)A(M,S)=E(M',P')E(M,S'), the AI=A(M,P)I(M,S) are Darii's premises, and the E'I'=E(M',P)I(M',S)=A(P,M)O(S,M) are Baroco's premises.

Then, one may use an algebraic method, (similar to the "set decompositions into subsets of U" à la Jevons, or similar to the trees, or, equivalently, similar to the subscripts' method, à la Carroll), to express all the 32 LC entailing PCPs and deduce their LCs via only four formulas. Another two methods for obtaining the LCs will be described in Sections 3 and 4.

Note that Christine Ladd, 1883, p. 26, lists, in her "incomplete exclusion operator" notation, V, and her "wholly exclusion operator" notation, $\overline{V}$, (identical with the present I, and resp., E, operators/statements), "eight propositions" which are in fact exactly the same statements, A,E,I,O, A',E',I',O', making the P and S cubes above. She cites DeMorgan, 1860, as using these "eight propositions", as well as for providing a definition of the affirmative and negative statements which depends on the number of the negative terms in an exclusion statement.

Each of the eight PCPs of type (1), (or type Barbara), E(M*,P*)E(M*',S*), entails two, non-independent, but different, LCs – since these PCPs are in fact "chain inclusion" sorites with just one middle term. By handling these PCPs using a tree like method, (as just noted above - very similar to Jevons' method of decompositions into subsets, or to Lewis Carroll's method of subscripts, and also similar to Carroll's own method of trees), and reading the sorite in the Aristotelian way, i.e. starting with the S* term of the PCP, one gets: S*=S*M*+S*M*'= S*M*= S*M*P*+ S*M*P*'= S*M*P*', where S*M*P*':= M*∩S*∩P*', etc., and the + sign denotes union of disjoint sets. Reading the sorite in the Goclenian way, i.e. starting with the P* term of the PCP, one gets the 2$^{nd}$ LC: P*=P*M*+ P*M*'= P*M*'= P*M*'S*+ P*M*'S*' =P*M*'S*'. Note that, based on the double inclusions S*⊆M*⊆P' and P*⊆M*'⊆S' – inferred from the premises E(M*,P*)E(M*',S*) - each LC implies the other, i.e., only one LC is independent, but the two LCs pinpoint to different subsets of U. Listing, after a column sign, and separated by semi-columns, all the possible LCs and LCs' formats, one may write:

**(α)** (Type 1) E(M*,P*)E(M*',S*): S*=S*M*P*', P*=P*M*'S*'; A(S*,P*') [=A(P*,S*')=E(S*,P*)] - after M* is dropped, (classical style!), the two LCs become identical and the number of all LCs of type (1) reduces to only four; I(S*, P*') - after ei on S*; I(P*, S*') - after ei on P*.

Thus the eight PCPs of type (1) generate 16 VCAs, out of which only 8 remain distinct after M* and M*' are dropped from their LCs; they also generate 16 ei VCAs – which remain all distinct – although there are only four distinct ei LCs out of which the middle term was dropped, since, as sets, I(S*, P*') and I(P*, S*'), are equal: {I(S*, P*'), P*∈{P, P'}, S*∈{S, S'}}={I(P*, S*'), P*∈{P, P'}, S*∈{S, S'}}.

It is even easier to find the LC of the eight type (2), or Darapti type, PCPs, E(M*,P*)E(M*,S*). They entail just one universal LC: M*=M*P* + M*P*'= M*P*'= M*P*'S* + M*P*'S*' =M*P*'S*', or A(M*,M*P*'S*') which reflects the inclusions M*⊆P*' and M*⊆S*' asserted by the premises. Thus, for the type (2) PCPs, the middle term itself is "the subject" of the precise LC; it can be eliminated only via ei on M*. Listing, after a column sign, and separated by semi-columns, all the possible LCs, one may write the 8 Darapti type (2) VCAs as:

**(β)** (Type 2) E(M*,P*)E(M*,S*): M*=M*P*'S*', A(M*,M*P*'S*'); M*P*'S*'≠Ø if M*≠Ø; I(S*',P*') if M*≠Ø and after M* is eliminated. The non-ei VCA whose universal LC is All M* is M*P*'S*', (resp. the ei VCA with the middle term M* eliminated from the LC), will be called the non-ei Darapti, (resp. the ei Darapti). The non-ei Darapti is not acknowledged by the CCS, because the middle term, M*, is the subject of the LC and thus can not be simply dropped from the LC. (Imagine a Barbari without a Barbara, or a Bramantip without - wait a minute - Bramanta? Now imagine an ei Darapti, without the non-ei Darapti. Perfect - welcome to CCS – the classical categorical syllogistic!) In conclusion, the eight PCPs of type (2) generate, before the middle term is dropped from the LCs, eight diferent non-ei VCAs and eight diferent ei VCAs. Once the middle term is dropped from the LCs, only eight ei VCAs remain.

There are eight distinct PCPs of type (3a), $E(M^*,P^*)$ $I(M^*,S^*)$, (or Darii/Datisi type), and eight distinct PCPs of type (3b), $I(M^*,P^*)$ $E(M^*,S^*)$, (or Disamis/Dimaris type). The very short trees "revealing" the LCs are: $\emptyset \neq M^*S^* = M^*S^*P^* + M^*S^*P^{*'} = M^*S^*P^{*'}$, and, resp., $\emptyset \neq M^*P^* = M^*S^*P^* + M^*S^{*'}P^* = M^*S^{*'}P^*$. One may also write:

**(γ)** (Type 3a)   $E(M^*,P^*)$ $I(M^*,S^*):M^*S^*P^{*'} \neq \emptyset$; $O(S^*,P^*)$ $[=I(S^*,P^{*'})]$ - after $M^*$ is dropped from the LC.

**(δ)** (Type 3b)   $I(M^*,P^*)$ $E(M^*,S^*):M^*S^{*'}P^* \neq \emptyset$; $O(P^*,S^*)$ $[=I(S^{*'},P^*)]$ - after $M^*$ is dropped from the LC.

Note that after ei is imposed on any one of the S,S',P,P', M,M' sets, left with only one possible non-empty subset (out of four) - as a consequence of the two universal premises of types either (1) or (2) having been taken into account, the resulting ei LC is identical with an LC of a PCP of type (3a) or (3b): one of the eight subsets of U is definitely not empty.

Types (1), (2), (3a) and (3b) PCPs, each contribute with eight PCPs to the 32 PCPs which entail at least one LC. If the precise, one subset pointing to, LCs are kept, then the total count of the VCAs generated by the previous four PCP types, comes to 40 non-ei VCAs and 24 ei VCAs. If the middle term is eliminated from the (precise) LCs then the total count of the VCAs generated by the previous four PCP types, shrinks to only 24 (non-ei) VCAs, while the total number of ei VCAs remains the same: 24. (As already noticed, the other 32 PCPs, of types (4a), (4b), (5a), (5b), containing eight PCPs each, do not entail any LC.)

One can see by inspection, that the above formulas (α), (β), (γ), (δ), restricted to LCs from which the middle term was eliminated, are in agreement with the four RofVCA listed below:

RofVCA #1: The distribution of the end terms is conserved in all non-ei and ei VCAs, except in type (1) ei VCAs, where ei on $S^*$, (resp. $P^*$) changes $S^*$, (resp. $P^*$), from distributed in the PCP to undistributed in the ei LC, while the distribution of the other end term, $P^*$, (resp. $S^*$), remains the same as it was in the PCP.

RofVCA #2: From two universal premises follows a universal LC, except when an ei condition on $M^*$, $S^*$ or $P^*$ is imposed – then a particular LC follows.

RofVCA #3: If the PCP contains one particular premise, then the LC is particular.

RofVCA #4: From two affirmative or two negative premises an affirmative LC follows; from one affirmative and one negative premises a negative LC results. This rule is valid even for LCs obtained after an ei condition was imposed.

Note that the formulas (α), (β), (γ), (δ) already prove, by inspection, the RofVCA #1, #2, and #3. But to prove the RofVCA #4 one has to show that RofVCA #4 is satisfied for each choice of $S^*,M^*,P^*$ in each of the formulas (α), (β), (γ), (δ). In the reverse direction, one can apply – to the PCPs of types (1), (2), (3a) and (3b) – the RofVCA listed above, to deduce, the portions of the formulas (α), (β), (γ), (δ) which contain only LCs from which the middle term was eliminated - see Sections 11 and 12.

Note that Aristotle's definition, (Striker 2009: p.20), "A syllogism is an argument in which, certain things being posited, something other than what was laid down results by necessity because these things are so", provides not only a characterization of a syllogism – both premises are necessary to obtain the LC and the LC has to validly follow from the premises - but also a justification, (or a pretext), for the elimination of the middle term from the LC. Nevertheless this elimination always weakens the LC, which instead of asserting something about a unique subset of U, will now assert the same thing, less precisely, about two subsets of U – namely that two subsets might be non-empty or that at least one of the two is definitely non-empty. Moreover, the contradictory statement of the weakened LC is stronger, (since it negates something about a larger number of sets), than the contradictory statement of the initial, stronger LC – which referred to a unique subset of U. (Example: compare "John lives in Miami" with "John lives on Earth": the negation of the less precise information sends John into outer space, while the negation of the stronger info about John, places him only out of Miami. It reflects the fact that negating a multiple Or statement produces a multiple And statement.) Thus, in CCS, when performing an indirect reduction, or a reductio ad absurdum proof of validity, one proves a weakened LC using stronger than necessary premises. For example, Darii's validity may be proved, by impossibility, from Camestres. But this is unnecessary: suppose, by impossibility, that Darii's precise LC, E(M,P') I(M,S):SPM≠Ø, is false, i.e., SPM=Ø (by the law of excluded middle). Then, from Darii's general premise, A(M,P), i.e., MP'=Ø, it results SM=SMP+ SMP' =Ø, which already contradicts Darii's particular premise, I(S,M) or SM≠Ø – no Camestres had to be invoked, and there is no need to suppose, (the stronger), SP=Ø, (the contradictory of Darii's weakened LC SP≠Ø), since supposing, by impossibility, that SPM=Ø, suffices.

### 3. STM and the Karnaugh map for n=3 (drawn for subsets instead of truth values)

The universal set U is graphed as a rectangle – but the left and right borders of the rectangle are glued together to generate a cylinder, so that S'PM and S'P'M, (resp. S'PM' and S'P'M'), are adjacent:

| S'P'M  | SP'M  | SPM  | S'PM  |
|--------|-------|------|-------|
| S'P'M' | SP'M' | SPM' | S'PM' |

*Fig. 1*

One may easily check that the eight subsets of the Marquand (1881) – Veitch (1952) - Karnaugh (1953) map for n=3, possess the same adjacencies as the seven subsets from the 3-circle Venn diagram, (Venn, 1880), to which one has added an eighth set, S'P'M' - drawn around the three S,P,M intersecting Venn circles. As mentioned, since the eight subsets of Figure 1 are the "elementary" subsets of U, one calls them just subsets; no other set will be a "subset". Note that for the STM's purposes, a rectangular Karnaugh map, (or K-map), suffices, since one deals with subsets not with rows of a truth value table "ordered in Grey code". One may consider that the n=3 Karnaugh map from Figure 1, was obtained by mirroring to the right, (i.e., by considering that the right edge of the map is a mirror), an n=2 Karnaugh map for S and M,

| S'M  | SM  |
|------|-----|
| S'M' | SM' |

and then indexing with P', (resp. P), the old, (resp. new) four cells.
Analogously an n=4 Karnaugh map will be obtained by mirroring toward the bottom of the page an n=3 Karnaugh map for S, M, and P, and then indexing, e.g., with $M_1$, (resp. $M_1'$), the old, (resp. new), 8 cells of the n=3 Karnaugh map:

| S'P'MM$_1$   | SP'MM$_1$   | SPMM$_1$   | S'PMM$_1$   |
|--------------|-------------|------------|-------------|
| S'P'M'M$_1$  | SP'M'M$_1$  | SPM'M$_1$  | S'PM'M$_1$  |
| S'P'M'M$_1$' | SP'M'M$_1$' | SPM'M$_1$' | S'PM'M$_1$' |
| S'P'MM$_1$'  | SP'MM$_1$'  | SPMM$_1$'  | S'PMM$_1$'  |

In fact it does not matter if one mirrors the n=3 Karnaugh map always toward the right of the page, thus making a Karnaugh map for any n have only two rows, (or even one row, if one always mirrors toward the right – see below), or one mirrors the n-1 Karnaugh map either toward the right or toward the bottom of the page – such that the Karnaugh maps for n=even are always square maps (before any "gluing" is done). After gluing together the left and right edges, (resp. the left and right edges, and top and bottom edges), of a Karnaugh map with only two rows, (resp. more than two rows), any two adjacent cells will differ by just one value of the index letters (or just one truth value in the string of truth values which identifies each cell – in the case of Karnaugh maps for truth values). For more details about Karnaugh maps on non-standard numbers of rows, see the Yasser-maps of Abdalla, 2015. A Yasser-map for n=3, on only one row, is shown below; it was started, for n=1, with

| S | S' |

and continued via mirroring towards the right two times, adding each time, in the leftmost position, one more index letter – not accented to the left of the mirroring edge, and accented to the right of the mirroring edge,

| MPS | MPS' | MP'S' | MP'S | M'P'S | M'P'S' | M'PS' | M'PS |

such that the letters indexing two adjacent cells differ in only one letter being switched from accented to non accented. (If each letter is replaced by the digit 0, and then 0' is replaced by 1, one gets a table with the cells ordered in Grey code from left to right: 0,1,3,2,6,7,5,4 – as in the Figure 2 of Abdalla, 2015.)

Representing Barbara's premises on such a map, one obtains again the LCs: S=SPM, P'=S'P'M'. The LC of any VCA will refer to only one cell of an n=3 Karnaugh map, and an LC of any of the sorites considered in Section 6 will refer to only one cell of a Karnaugh map for an n equal to the number of terms in the sorite.

On the Figure 1 map above, after Barbara's premises empty 4 subsets out of 8, the other 4 subsets left, are still connected on a Karnaugh cylindrical map, but on a rectangular map, S'P'M'=P' is disconnected from S, M and P, which still satisfy S$\subseteq$ M$\subseteq$ P; M' is made of two disconnected subsets, P' and S'M'P.) On the Figure 1 map above, one can see that any universal premise empties two "horizontal subsets" located either on the M or the M' row, and its contradictory particular premise places set elements in at least one of the same two "horizontal subsets" emptied by the contradictory universal premise. This makes it clear that a PCP made of two particular premises does not entail any LC (beyond a possible verbatim repetition of premises' content). Also, it makes it clear that there is no LC entailed by any PCP made of one universal premise emptying two subsets on M', and a particular premise "acting on M", i.e. placing set elements into two subsets of M, (or the other way around). The Karnaugh map also shows that in STM, (and thus in CCS, too), any LC is characterized by the fact that it singles out just one of the eight subsets partitioning U, and that only four types of LCs are possible:

(1) One of the sets S or S' and one of the sets P or P', have each three empty subsets and one possibly not empty subset, (2) One of the sets M or M' has three empty subsets and one possibly not empty subset, (3a) and (3b) One of the subsets of U is definitely not empty. The above LC types are entailed by the already mentioned four PCP types: (1) PCPs made of two universal premises, that empty a total of four subsets on different rows - M and M'. In this case the intersecting set of the two sets having three empty subsets and one possibly not empty subset is itself empty – and this is the LC retained by the VS Barbara, Celarent/Cesare and Camestres/Camenes. When, e.g., Celarent's premises empty three subsets of S, they also empty three subsets of P – hence the two LCs entailed by only one PCP; in turn, each LC implies that S∩P is empty. (2) PCPs made of two universal premises, that empty a total of three subsets on the same row – either M or M'. (3a) and (3b) PCPs made of one universal plus one particular premises, one emptying, and one placing elements into, subsets of the same row.

One can notice that, e.g., a (3a) PCP, E(M,P')I(M,S), Darii/Datisi, (obtained when choosing S*=S, M*=M, P*=P' in (3a)), and a (3b) PCP, I(M,P)E(M,S'), Disamis/Dimaris, (choose S*=S', M*=M,

P*=P in (3b)), have the same LC: SMP≠Ø. Analogously Bocardo, O(M,P)A(M,S)=I(M,P')E(M,S') and Ferio/Festino/Ferison/Fresison, E(M,P)I(M,S), have the same LC: SMP'≠Ø. Note again that all the (3a) VCAs become the (3b) VCAs via the relettering S↔P, i.e., on specific VCA examples, to keep up with the convention of listing firstly the P- premise, one just have to switch the premises' order to obtain a (3b) VCA from a type (3a) VCA, and vice versa. By working with 64 distinct PCPs, STM "raises to an equal footing" any of the eight subsets partitioning U. Because of this, each one of the eight subsets of U, S*P*M*, appears as the "subject" of an LC exactly eight times: three times as a "possibly non empty subset" - in LCs of types (1) and (2), S*=S*P*M*, P*=S*P*M*, (these are LCs of the eight PCPs of type (1), each entailing two LCs per PCP), and M*=S*P*M*, (this is the LC of one of the eight type (2) PCP), and five times as a " definitely non empty subset" LC, (twice as the LC S*P*M*≠Ø of PCPs of types (3a) and (3b), and 3 times as the ei LC generated by the previous three "universal LCs" using ei on S*, P* and resp. M*. This shows that STM "likes equally" each of the U's eight subsets. By contrast, CCS "likes" only three subsets, SPM, SP'M, and SP'M': any VS LC asserts that either one of the above three subsets is not empty, (if the PCP contains one particular and one universal premises, both acting on the same set - either M or M'), or that one of the above three subsets is the only part of either S,P or M, which is possibly left non empty by a PCP containing two universal premises. More precisely, SPM appears as LC of Barbara, Barbari, Bramantip, (but "Bramanta", having A(P,S) as LC is not a VS, since its LC differs from the standard A(S,P)), Darapti, Darii/Datisi, Disamis/Dimaris; SP'M appears as LC of Celarent/Cesare, Celaront/Cesaro, Felapton/ Fesapo, Bocardo, Ferio/Festino/Ferison/Fresison; SP'M' appears as LC of Camestres/Camenes, Camestros/Camenos and Baroco.

## 4. The Euler diagrams method for finding the LCs

Another "LC finding method" is based on a specific choice of Euler diagrams: as in the previous two methods, any particular premise is represented as a non-empty set intersection, while any universal premise is represented not as an empty set intersection, but as a set inclusion of adeptly chosen sets, such that all PCPs from each of the four VCA types comport with the same type of graph. The usual 3-circle Venn diagram representation of U, see Quine, pp. 98-108, or Hurley, pp. 272-275, and the usual mood and figure CCS classification of the VS, do not illuminate the fact that the VS are of four types, too, and that the VS of the same type are equivalent. Graphing a PCP on a 3-circle Venn diagram does not offer any LC insight – except that if the 3-circle graph of the PCP validates the proposed LC, then the LC is considered (graphically) proven. Using Euler diagrams, e.g., O(P,M):= "Some P is not M" is represented as two intersecting circles corresponding to P∩M'≠Ø; this also shows that even if O(P,M) is formulable using only the positive terms P and M, in fact O(P,M) "acts" on M', since it lays set elements in the intersection PM':=P∩M'. One represents A(M,P):="All M is P" as an M-circle included in the P-circle, and A(P,M) as either a P-circle included in the M-circle, or, equivalently, by contraposition, as an M'-circle included in the P'-circle, (in accord also to the fact that A(P,M) means PM':=P∩M'=Ø, i.e., A(P,M), acts on M' by emptying the PM' intersection; moreover P' is the largest part of U which does not intersect P – thus M'⊆P' since PM'=Ø. The E(M,P) statement is represented in the previous two methods as MP=Ø, while in the Euler diagram method is represented as either an M-circle included in the P'-circle, or, as a P-circle included in the M'-circle. One easily finds that the graphical patterns generated by the Euler diagram method are as follows: Any type (1) PCP can be graphed as a double inclusion with M, (or M'), in the middle, α⊆M*⊆β, where α and β are two different letters chosen from the possible "end terms" {S,S',P,P'}, where M*∈{M, M'}, and the LC can be written as an equality α=α∩M*∩β=:αM*β, or, All α is M*β. After the middle term M* is eliminated/dropped, this may be written All α is β. (Replacing α⊆M*⊆β by the equivalent inclusions β'⊆M*'⊆α', one gets a second LC: β'=α'∩M*'∩β'=:α'M*'β'.) Any type (2) PCP can be graphed as M*⊆α, M*⊆β, where the LC, again, can be written as an equality M*= M*αβ or All M* is α∩β=:αβ. After the ei condition, M*≠∅, is imposed, the LC may be written as I(α,β). Any type (3a) or (3b) PCP can be graphed as M*∩α=:αM*≠∅, M*⊆β, where the LC can now be expressed as an "inequality" αM*β≠∅, or, as Some α is β - after the middle term M* is, again, dropped from the LC. Note that the precise LC always contains M or M', and that the inclusions characteristic to each VCA type also give the recipe for finding the LCs: for types (1) and (2) write the LC as an equality of the set appearing on the left side of the inclusions, to, on the right hand side of the equality, the intersection of all the sets appearing on the right hand side of the inclusions with the set appearing on the left side of the inclusions, thus giving the LCs, α=αM*β, and resp., M*= M*αβ. The LC of a type (3a) or (3b) VCA starts with the intersection of the middle term M*, (either M or M'), and an end term being non-empty; the inclusion of M* in the other end term, adds this latter end term to the initial 2-set intersection (which already contains M*). All four VCA types can be generalized to sorites, (see Section 6), whereby the corresponding inclusions and LCs will be:

Type (1) α⊆$M_1$⊆$M_2$⊆$M_3$⊆...⊆$M_n$⊆β; with the LC α=α$M_1M_2M_3...M_n$β. Notice that the first chain of inclusions is equivalent with this second chain: β'⊆$M_n$'⊆$M_{n-1}$'⊆$M_{n-2}$'⊆...⊆$M_1$'⊆α', whose LC is β'=α'$M_1$'$M_2$'$M_3$'...$M_n$'β'. Obviously, α≠β', but A(α,β)=A(β',α') – after all the "middle terms" are dropped.

Type (2) M⊆α, M⊆$M_1$, M⊆$M_2$, M⊆$M_3$, ..., M⊆$M_n$, M⊆β; with the LC M=Mα$M_1M_2M_3...M_n$β

Type (3a) and (3b) M∩α=:Mα≠∅, M⊆$M_1$, M⊆$M_2$, M⊆$M_3$, ..., M⊆$M_n$, α⊆$α_1$, α⊆$α_2$, α⊆$α_3$, ..., α⊆$α_r$ ; with the LC Mα$M_1M_2M_3...M_nα_1α_2α_3$, ...,$α_r$ ≠∅. Sorites with (slightly) more complex structures are presented in Section 6.

## 5. The "structure" of the 6 by 6 matrix of the PCPs formulable only via positive terms

One has five PCPs of type (4a), made of two particular premises both acting on either M or M': I(M,P*) I(M,S*), P*∈{P, P'}, S*∈{S, S'}, and also I(M',P)I(M',S) – since there exists only one PCP formulable via positive terms if both premises contain M'. There are only four PCPs of type (4b): I(M,P*) I(M',S), P*∈{P, P'}, and I(M',P) I(M,S*), S*∈{S, S'}. There are four PCPs of type (5a), E(M,P*) I(M',S), P*∈{P, P'}, and E(M',P) I(M,S*), S*∈{S, S'}, and another four PCPs of type (5b): I(M,P*) E(M',S), P*∈{P, P'}, and I(M',P) E(M,S*), S*∈{S, S'}. Each of the remaining 19 PCPs entails at least one LC. These PCPs generate VS if the LCs are one of the A(S,P), E(S,P), O(S,P), I(S,P) statements, or generate other VCAs which CCS discards for reasons mentioned below. Firstly the PCP, E'E'=E(M',P)E(M',S)= A(P,M)A(S,M), is discarded because M is undistributed in both premises; then, the PCPs, EE, OE, EO, are discarded because in each PCP both premises are negative. (Here is the RofVS #1: The middle term has to be distributed in at least one premise; and the RofVS #2: two negative premises are not allowed. The above four PCPs generate VCAs or ei VCAs, all having I(S',P') as the LC. For a complete discussion, please see Section 11.) The PCPs, AI (Darii), IA (Disamis), AA (Darapti), have I(S,P) as the LC, or the ei LC (in Darapti's case), and EE' (Celarent), and E'E (Camestres) have E(S,P) as the LC. There are five PCPs which have either A(S,P) as LC, AE' (Barbara), or have O(S,P) as LC: EI (Ferio), OA (Bocardo), E'I' (Baroco), EA (Felapton). The last five PCPs out of the 19 which entail LCs, are dismissed by CCS since their LCs are, A(P,S) for E'A (Bramanta?! - which by ei on P gives the acceptable ei LC I(S,P) – Bramantip), and O(P,S) for IE (Fireo?), AO (Bacordo?), I'E' (Boraco?), AE (Falepton?). The CCS' reason for dismissal of the VCAs having A(P,S) and O(P,S) as LCs is that a premises' transposition and a relettering S↔P would transform them into VS. Another very good reason to dismiss the VCAs whose LCs are A(P,S) and O(P,S), is that their premises, "contradict pairwise" the premises of the VS whose LCs are A(S,P) and O(S,P). Indeed, the following pairs of PCPs have premises which either contradict each other, EI and IE, OA and AO, E'I' and I'E', or are contraries to each other, EA and AE – the latter two PCPs are both sound only if M=P=S=0. Barbara's and Bramanta's premises can be all sound only if S=P=M. Also note that the dismissed PCPs, IE, E'A, AE, AO, I'E', and their LCs, satisfy all the RofVS.

One could have classified the 19 PCPs entailing at least one LC according to their types. Type (1) PCPs are Barbara, Bramanta, Celarent and Camestres. One can see that any of the six PCP pairs one may consider, imposes, if all the four premises of such a pair are true, a particular structure on the the universe of discourse, U. For example, if the premises of Celarent and Camestres are all true, i.e., MP=0, SM'=0, and also, PM'=0, SM=0, it will result that P=S=0. If Barbara and Celarent have both true premises, i.e., MP'=0, SM'=0, and also, MP=0, SM'=0, it results M=MP'+MP=0; therefore M'=U, and SM'=SU=S=0. This shows that, once the middle term is specified, only one PCP of type (1) can have sound premises without a particular structure being imposed on U. The PCPs of type (2) are: AA (Darapti), EA (Felapton), AE (Falepton), EE, and E'E'. Again, one can see that any of the ten PCP pairs one may consider, imposes, if all the four premises of such a pair of PCPs are true, a particular structure on the the universe of discourse, U. For example, if the premises of Darapti and Felapton are all true, i.e., MP'=0, MS'=0, and also, MP=0, MS'=0, it will result that M=0, and the ei condition, M≠Ø, can not be imposed. If the premises of Darapti and EE are all four true, MP'=0, MS'=0, MP=0, MS=0, it follows again that M=0. One can again show that, unless U has a particular structure, only one PCP of type (2) can have sound premises – once the middle term is specified. The above arguments are based on the fact that any two of the four universal P-premises, (resp. S-premises), being simultaneously true already imposes a particular structure on the universal set U. (For example E(M,P) and A(M,P) imply M=Ø; E(M,P) and A(M',P) imply M=P', etc.) See Section 13 for complete details. The PCPs of type (3a) are: AI (Darii), EI (Ferio), E'I' (Baroco), AO (Bacordo?), EO. The PCPs of type (3b) are: IA (Disamis), IE (Fireo?), I'E' (Boraco?), OA (Bocardo), OE.

## 6. Biliteral Sorites

The four types of PCPs which generate VCAs may be immediately extended to simple sorites of the Barbara, Darapti, Darii, and Disamis types. Since Darii and Disamis types are so similar, one discusses only the Darii type sorites:

1. The premises of a simple Barbara type Aristotelian/Goclenian sorite are,

SOR1:=$E(S^*, M_1^{*'})E(M_1^*{}', M_2^*)E(M_2^*{}', M_3^*)...E(M_i^*{}', M_{i+1}^*)...E(M_k^*{}', M_{k+1}^*)...E(M_n^*{}', P^*)=$
    $A(S^*, M_1^*{}')A(M_1^*{}', M_2^*{}') A(M_2^*{}', M_3^*{}')...A(M_n^*{}', P^*{}')$, with:

LC$_1$: $S^*=S^*M_1^*{}'M_2^*{}'M_3^*{}'...M_i^*{}'...M_k^*{}'...M_n^*{}'P^*{}'$.  (the Aristotelian LC)

LC$_2$: $P^*=P^*M_n^*...M_k^*...M_i^*...M_3^*M_2^*M_1^*S^*{}'=P^*M_1^*M_2^*M_3^*...M_i^*...M_k^*...M_n^*S^*{}'$ (the Goclenian LC of the same premises. As in Barbara's PCP which contains just two chain links, E(S,M') and E(M,P), the two LCs are not independent because they are based on two chains of inclusions which are deducible one from the other: $P^*\subseteq M_n^*{}'\subseteq...\subseteq M_3^*\subseteq M_2^* \subseteq M_1^*\subseteq S^*{}'$ and $S^*\subseteq M_1^*{}'\subseteq M_2^*{}' \subseteq M_3^*{}'\subseteq...\subseteq M_n^*{}'\subseteq P^*{}'$.)

2. The premises of a simple Darapti type sorite are,

SOR2:=$E(S^*, M^*)E(M_1^*, M^*)E(M_2^*, M^*)...E(M_i^*, M^*)E(M_{i+1}^*, M^*)...E(M_n^*, M^*) E(M^*, P^*)$, with the LC, $M^*=M^*S^*{}'M_1^*{}'M_2^*{}'M_3^*{}'...M_i^*{}'M_{i+1}^*{}'...M_n^*{}'P^*{}'$, where, if $M^*\neq\emptyset$, some of the many ei LCs are $MS^*{}'M_1^*{}'M_2^*{}'M_3^*{}'...M_i^*{}'M_{i+1}^*{}'...M_n^*{}'P^*{}' \neq\emptyset$, I(S*,P*'), etc.

3. The premises of a simple Darii type sorite are:

$SOR3:=E(S^*,M_x^*)E(M_1, M_y^*)E(M_2^*, M_y^*)...I(M_x^*, M_y^*)...E(M_i^*, M_y^*)E(M_{i+1}^*,M_x^*)... E(M_n^*, M_y^*)$
$E(M_x^*, P^*)$

LC: $M_x^*M_y^*S^{*\prime}M_1^{*\prime}M_2^{*\prime}M_3^{*\prime}...M_i^{*\prime}M_{i+1}^{*\prime}...M_n^{*\prime}P^{*\prime} \neq \emptyset$

Note that there are two types of "chains" making up the three sorites above. One is a D-chain, (or Darapti-chain), $E(\alpha,\beta)E(\alpha,\gamma)E(\alpha,\delta)$, the other is a B-chain, (or Barbara-chain), $E(a,b)E(b',c)E(c',d)$, where $\alpha,\beta,\gamma,\delta,a,b,c,d$ are sets, $\alpha\cap\beta\cap\gamma\cap\delta$ is written $\alpha\beta\gamma\delta$, etc., and the LCs are $\alpha=\alpha\beta'\gamma'\delta'$ and, respectively, $a=ab'c'd'$ and $d=dcba'$.

Note that in SOR3 one could have had just one D-chain instead of two, (for example, replace $M_y^*$, everywhere but in $I(M_x^*, M_y^*)$, with $M_x^*$ - to obtain just one D-chain), or one might have had one or two B-chains, (or one B-chain and one D-chain), originating at $M_x^*$ or/and $M_y^*$.

To the two 3-links B and D-chains just above, one may add many other premises/chains and still have the extended premises entailing at least one LC:

$E(a,b)E(b',c)E(c',d)$ $E(a,b_1)E(b_1',c_1)E(c_1',d_1)...E(a,b_n)E(b_n',c_n)E(c_n',d_n)$, entails the LC $a=ab'c'd'b_1'c_1'd_1'...b_n'c_n'd_n'$, after all $a=ab_k'c_k'd_k'$, $k=1,2,...,n$, are substituted in (the initial) $a=ab'c'd'$

If in the above premises one replaces all $d_i$ $i=1,2,...,n$, by d, then the LC becomes $a=ab'c'd'b_1'c_1'...b_n'c_n'$, and there is a second LC: $d=da'bcb_1c_1d_1...b_nc_n$.

Another "extended premises" possibility, would have been to start the B-chains not at the term a, but at any other term, e.g., c':

$E(a,b)E(b',c)E(c',d)$ $E(c',b_1)E(b_1',c_1)E(c_1',d_1)$ $E(c',b_2)E(b_2',c_2)E(c_2',d_2)...E(c',b_n)E(b_n',c_n)E(c_n',d_n)$, entails the LC $a=ab'c'd'b_1'c_1'd_1'...b_n'c_n'd_n'$, obtained by substituting the c' term in $a=ab'c'd'$ with $c'=c'b_1'c_1'd_1'...b_n'c_n'd_n'$. If one would have started some B-chains at c', and some other B-chains, (or a D-chain), at c, then no LC pointing to the unique subset of the entire universe of discourse, U, would have been possible. A "substitution principle" is necessary to hold if it were that a set of premises entails an LC: the partial LCs of all the partial B and D-chains have to allow term substitutions until all partial LCs "unite" into an LC pinpointing to just one subset of the entire universe of discourse U comprising all the terms appearing in all the premises. All the subsets of U can be represented on a rectangular Karnaugh map drawn for subsets, not truth values (for which a toroidal K-map is preferable).

Other examples (involving more indexes):

One may say that each of the three (or four) type sorites admits "Darapti or/and Barbara decorations", i.e., a D-chain or/and B-chain sequence of premises may be added using any term/letter as the "starting base" - exactly as the term (or letter) M was used as the "base" of SOR2 - the Darapti type (2) sorite above.

For example, if to SOR1 one adds Darapti decorations of length=1 on $S^*$, $M_1^{*\prime}$, and $P^{*\prime}$, respectively, and a Darapti decoration of length=3 on $M_2^{*\prime}$, i.e., if one adds these premises: $E(S^*, M_{n+1}^*)E(M_1^{*\prime}, M_{n+2}^*)E(M_2^{*\prime}, M_{n+3}^*)E(M_2^{*\prime}, M_{n+4}^*)E(M_2^{*\prime}, M_{n+5}^*)E(P^{*\prime}, M_{n+6}^*)$, then one obtains a sorite which has just one LC not two. This LC is an "extension" of $LC_1$: $S^*=S^*M_1^{*\prime}M_2^{*\prime}M_3^{*\prime}...M_n^{*\prime}P^{*\prime}M_{n+1}^{*\prime}M_{n+2}^{*\prime}M_{n+3}^{*\prime}M_{n+4}^{*\prime}M_{n+5}^{*\prime}M_{n+6}^{*\prime}$.

For the Goclenian reading of SOR1 one needs to add, as Darapti decorations of the same lengths as above, e.g., these premises:

$E(S^{*\prime}, M_{n+1}^*)E(M_1^*, M_{n+2}^*)E(M_2^*, M_{n+3}^*)E(M_2^*, M_{n+4}^*)E(M_2^*, M_{n+5}^*)E(P^*, M_{n+6}^*)$, and then, again, one obtains a sorite which has just one LC not two. This LC is an "extension" of $LC_2$: $P^*=P^*M_1^*M_2^*M_3^*...M_n^*S^{*\prime}M_{n+1}^{*\prime}M_{n+2}^{*\prime}M_{n+3}^{*\prime}M_{n+4}^{*\prime}M_{n+5}^{*\prime}M_{n+6}^{*\prime}$.

One may add Darapti decorations to a Darapti type sorite, SOR2, via adding, e.g., these premises:
$E(S^{*\prime}, M_{n+1}^*)E(M_1^*, M_{n+2}^*)E(M_2^*, M_{n+3}^*)E(M_2^*, M_{n+4}^*)E(M_2^*, M_{n+5}^*)E(P^*, M_{n+6}^*)$.

Then one obtains a sorite whose LC is an "extension" of the above SOR2 LC:

$M^* = M^*S^{*'}M_1^{*'}M_2^{*'}M_3^{*'}...M_i^{*'}M_{i+1}^{*'}...M_n^{*'}P^{*'}M_{n+1}^{*'}M_{n+2}^{*'}M_{n+3}^{*'}M_{n+4}^{*'}M_{n+5}^{*'}M_{n+6}^{*'}$. Analogously, Darapti decorations of the same lengths as above may be added to a Darii type sorite, SOR3, via adding these premises:

$E(S^{*'}, M_{n+1}^*)E(M_1^{*'}, M_{n+2}^*)E(M_2^{*'}, M_{n+3}^*)E(M_2^{*'}, M_{n+4}^*)E(M_2^{*'}, M_{n+5}^*)E(P^{*'}, M_{n+6}^*)$.

Then one obtains a sorite whose LC is an "extension" of the SOR3 LC: $M_x * M_y *S^{*'}M_1^{*'}M_2^{*'}M_3^{*'}...M_i^{*'}M_{i+1}^{*'}...M_n^{*'}P^{*'}M_{n+1}^{*'}M_{n+2}^{*'}M_{n+3}^{*'}M_{n+4}^{*'}M_{n+5}^{*'} M_{n+6}^{*'}$.

One thus sees that to any term/letter appearing in the LC of a sorite of types 1, 2 or 3, one may add a "Darapti decoration", i.e., a (no matter how long) Darapti type sequence of premises with that letter as the "base", exactly as the term (or letter) M was the "base" of the Darapti type SOR2. "Barbara decorations" can be added in a similar way.

If represented on a subsets' Karnaugh map with enough variables, the LCs of the above SOR1, (resp. SOR2), say that $S^*$, $P^*$, (resp. $M^*$), were reduced by the universal premises to just one, possibly not empty, partition subset of the universal set U. The LC of SOR3 affirms that there exists a partition subset of the universal set U which is not empty.

A product shorthand notation of SOR1 would have as factors the chain links $E(M_i^{*'}, M_{i+1}^*)$, where i runs from 0 to n, and where $M_0^{*'} = S$, $M_{n+1}^* = P$. In this notation one can see that a B-chain can split, at any link i, into any number of new chains, originating either at $M_i^*$ or at $M_{i+1}$. For example, one may add to the SOR1 above, either premises which "come back in a loop" to the initial B-chain, such as in (a), $E(M_i^{*'}, \alpha_1)E(\alpha_1', \alpha_2)E(\alpha_2', \alpha_3)E(\alpha_3', M_k^*)$, or continue with other partitioning sets $\alpha_4$, $\alpha_5$, etc., from the universal set U, as in (b), $E(M_i^{*'}, \alpha_1)E(\alpha_1', \alpha_2)E(\alpha_2', \alpha_3)E(\alpha_3', \alpha_4) E(\alpha_4', \alpha_5)E(\alpha_5', \alpha_6)$. Or, new B-chains can originate at $M_{i+1}^*$: (c) $E(M_{i+1}^*, \beta_1)E(\beta_1', \beta_2)E(\beta_2', \beta_3) E(\beta_3', M_k^{*'})$ - these are premises which "come back in a loop", or/and, (d), $E(M_{i+1}^*, \gamma_1)E(\gamma_1', \gamma_2) E(\gamma_2', \gamma_3)E(\gamma_3', \gamma_4)E(\gamma_4', \gamma_5)E(\gamma_5', \gamma_6)$, or/and, (e), $E(M_{i+1}^*, \delta_1)E(\delta_1', \delta_2) E(\delta_2', \delta_3)E(\delta_3', \delta_4)E(\delta_4', \delta_5)$. As before, the "extended" SOR1a, SOR1b, SOR1c, SOR1d, and SOR1e, entail an LC only if, (as a result of all the new and old (universal) premises being applied), one of the n terms which generate a partition of U into $2^n$ subsets, has $2^n-1$ empty subsets and only one possibly non-empty subset of the new, $2^n$ "extended" partitioning of U, determined by all the terms, old and new, contained in all the premises, old and new - or "extended". For example, the extended sorite,

$SOR1a = E(S^*, M_1^*)E(M_1^{*'}, M_2^*)E(M_2^{*'}, M_3^*)...E(M_i^{*'}, M_{i+1}^*)...E(M_k^{*'}, M_{k+1}^*)...E(M_n^{*'}, P^*)$
$E(M_i^{*'}, \alpha_1)E(\alpha_1', \alpha_2)E(\alpha_2', \alpha_3)E(\alpha_3', M_k^*)$, still entails two LCs,

 $LC_{1a}$: $S^* = S^*M_1^{*'}M_2^{*'}M_3^{*'}...M_i^{*'}...M_k^{*'}...M_n^{*'}P^{*'}\alpha_1'\alpha_2'\alpha_3'$ (the Aristotelian LC obtained after $M_i^{*'} = M_i^{*'}\alpha_1'\alpha_2'\alpha_3'M_k^{*'}$ was substituted in the $LC_1$ of SOR1).

 $LC_{2a}$: $P^* = P^*M_n^*...M_k^*...M_i^*...M_3^*M_2^*M_1^*S^{*'}\alpha_1\alpha_2\alpha_3$ (the Goclenian LC obtained after $M_k^* = M_k^*\alpha_1\alpha_2\alpha_3M_i^*$ was substituted in the $LC_2$ of SOR1). The same is true for SOR1c – it will still admit two LCs: one has to substitute $M_{i+1}^*$ in $LC_2$, and $M_k^{*'}$ in $LC_1$ – see below. But note that if the last link of the "extension" (a) would have been $E(\alpha_3', M_k^{*'})$ instead of $E(\alpha_3', M_k^*)$, then the substitution $M_i^{*'} = M_i^{*'}\alpha_1'\alpha_2'\alpha_3'M_k^*$ in $LC_{1a}$ would have produced the empty set: $S^* = S^*M_1^{*'}M_2^{*'}M_3^{*'}...M_i^{*'}...M_k^*M_k^{*'}...M_n^{*'}P^{*'}\alpha_1'\alpha_2'\alpha_3'$. And no substitution would have been available to the Goclenian solution, $LC_2$ of SOR1, to account for the extension $M_k^{*'} = M_k^{*'}\alpha_1\alpha_2\alpha_3M_i^*$, while substituting $M_k^{*'}$ into $LC_1$ would have resulted again in an empty set.

The extended sorite,
$SOR1b = E(S^*, M_1^*)E(M_1^{*'}, M_2^*)E(M_2^{*'}, M_3^*)...E(M_i^{*'}, M_{i+1}^*)...E(M_k^{*'}, M_{k+1}^*)...E(M_n^{*'}, P^*)$
$E(M_i^{*'}, \alpha_1)E(\alpha_1', \alpha_2)E(\alpha_2', \alpha_3)E(\alpha_3', \alpha_4) E(\alpha_4', \alpha_5)E(\alpha_5', \alpha_6)$, entails only one LC, $LC_{1b}$:
$S^* = S^*M_1^{*'}M_2^{*'}M_3^{*'}...M_i^{*'}...M_k^{*'}...M_n^{*'}P^{*'}\alpha_1'\alpha_2'\alpha_3'\alpha_4'\alpha_5'\alpha_6'$ (the Aristotelian LC), while the Goclenian LC of SOL1, $LC_2$: $P^* = P^*M_n^*...M_k^*...M_i^*...M_3^*M_2^*M_1^*S^{*'}$, together with $\alpha_6 = \alpha_6\alpha_5\alpha_4\alpha_3\alpha_2\alpha_1M_i^*$ do not pinpoint to a unique subset of U. A Goclenian $LC_{2b}$ does not exist.

As already mentioned, the extended sorite SOR1c, similar to SOR1a, entails two LCs,

SOR1c=$E(S^*, M_1^*)E(M_1^{*\prime}, M_2^*)E(M_2^{*\prime}, M_3^*)...E(M_i^{*\prime}, M_{i+1}^*)...E(M_k^{*\prime}, M_{k+1}^*)...E(M_n^{*\prime}, P^*)$
$E(M_{i+1}^*,\beta_1)E(\beta_1',\beta_2)E(\beta_2',\beta_3)E(\beta_3',M_k^{*\prime})$,

    $LC_{1c}$: $S^*=S^*M_1^{*\prime}M_2^{*\prime}M_3^{*\prime}...M_i^{*\prime}...M_k^{*\prime}...M_n^{*\prime}P^{*\prime}\beta_1\beta_2\beta_3$ (the Aristotelian LC, obtained after substituting $M_k^{*\prime}=M_k^{*\prime}\beta_1\beta_2\beta_3 M_{i+1}^{*\prime}$ in $LC_1$).

    $LC_{2c}$: $P^*=P^*M_n^*...M_k^*...M_i^*...M_3^*M_2^*M_1^*S^{*\prime}\beta_1'\beta_2'\beta_3'$ (the Goclenian LC, obtained after substituting $M_{i+1}^*=M_{i+1}^*\beta_1'\beta_2'\beta_3'M_k^*$ in $LC_1$).

The extended sorite,
SOR1bd=$E(S^*, M_1^*)E(M_1^{*\prime}, M_2^*)E(M_2^{*\prime}, M_3^*)...E(M_i^{*\prime}, M_{i+1}^*)...E(M_k^{*\prime}, M_{k+1}^*)...E(M_n^{*\prime}, P^*)$
$E(M_i^*,\alpha_1)E(\alpha_1',\alpha_2)E(\alpha_2',\alpha_3)E(\alpha_3',\alpha_4) E(\alpha_4',\alpha_5)E(\alpha_5',\alpha_6) E(M_{i+1}^*,\gamma_1)E(\gamma_1',\gamma_2) E(\gamma_2',\gamma_3)E(\gamma_3',\gamma_4)E(\gamma_4',\gamma_5)$,
does not have any LC, since the two "LC partial extensions" $M_i^{*\prime}=M_i^{*\prime}\alpha_1'\alpha_2'\alpha_3'\alpha_4'\alpha_5'\alpha_6'$ and $M_{i+1}^{*\prime}=M_{i+1}^{*\prime}\gamma_1'\gamma_2'\gamma_3'\gamma_4'\gamma_5'$ can not be simultaneously substituted in either $LC_1$ or $LC_2$ of SOR1.
Finally, the extended sorite,
SOR1de=$E(S^*, M_1^*)E(M_1^{*\prime}, M_2^*)E(M_2^{*\prime}, M_3^*)...E(M_i^{*\prime}, M_{i+1}^*)...E(M_k^{*\prime}, M_{k+1}^*)...E(M_n^{*\prime}, P^*)$
$E(M_{i+1}^*,\gamma_1)E(\gamma_1',\gamma_2) E(\gamma_2',\gamma_3)E(\gamma_3',\gamma_4)E(\gamma_4',\gamma_5)E(M_{i+1}^*,\delta_1)E(\delta_1',\delta_2) E(\delta_2',\delta_3)E(\delta_3',\delta_4)E(\delta_4',\delta_5)$, has one LC only - $LC_{1de}$ does not exist since neither of the "extended partial" LCs, $M_{i+1}^{*\prime}=M_{i+1}^{*\prime}\gamma_1'\gamma_2'\gamma_3'\gamma_4'\gamma_5'$ and $M_{i+1}^{*\prime}=M_{i+1}^*\delta_1'\delta_2'\delta_3'\delta_4'\delta_5'$ can be substituted in the $LC_1$, $S^*=S^*M_1^{*\prime}M_2^{*\prime}M_3^{*\prime}...M_i^{*\prime}...M_k^{*\prime}...M_n^{*\prime}P^{*\prime}$, of SOR1.
But, via substitutions, one obtains $LC_{2de}$: $P^*=P^*M_n^*...M_k^*...M_i^*...M_3^*M_2^*M_1^*S^{*\prime}\gamma_1'\gamma_2'\gamma_3'\gamma_4'\gamma_5'\delta_1'\delta_2'\delta_3'\delta_4'\delta_5'$, (the Goclenian LC of SOR1).
Anywhere in the D-chain of SOR2 one may have added an arbitrary number of B-chains, such as the above: (a) $E(M_i^*,\alpha_1)E(\alpha_1',\alpha_2)E(\alpha_2',\alpha_3)E(\alpha_3',M_k^*)$, or/and (b), $E(M_i^*,\alpha_1)E(\alpha_1',\alpha_2)E(\alpha_2',\alpha_3)E(\alpha_3',\alpha_4)$ $E(\alpha_4',\alpha_5)E(\alpha_5',\alpha_6)$, or/and (c') $E(M_{i+1}^{*\prime},\beta_1)E(\beta_1',\beta_2)E(\beta_2',\beta_3) E(\beta_3',M_k^*)$, or/and, (d'), $E(M_{i+1}^{*\prime},\gamma_1)$ $E(\gamma_1',\gamma_2)E(\gamma_2',\gamma_3)E(\gamma_3',\gamma_4) E(\gamma_4',\gamma_5)$, or/and, (e'), $E(M_{i+1}^{*\prime},\delta_1)E(\delta_1',\delta_2) E(\delta_2',\delta_3)E(\delta_3',\delta_4)E(\delta_4',\delta_5)$, since all the "extended partial" LCs, $M_i^{*\prime}=M_i^{*\prime}\alpha_1'\alpha_2'\alpha_3'M_k^{*\prime}$, $M_i^{*\prime}=M_i^{*\prime}\alpha_1'\alpha_2'\alpha_3'\alpha_4'\alpha_5'\alpha_6'$, $M_{i+1}^{*\prime}=M_{i+1}^{*\prime}\gamma_1'\gamma_2'\gamma_3'\gamma_4'\gamma_5'$, and $M_{i+1}^{*\prime}=M_{i+1}^{*\prime}\delta_1'\delta_2'\delta_3'\delta_4'\delta_5'$ can be substituted in the LC of SOR2.

## 7. Inside each of the four types of VCA, any VCA can be re-written, without changing its PCP's nor its LC's "contents", as any other VCA of the same type

The short way to show that each of the four VCA types can be written, respectively as a Barbara, Darapti, Darii and Disamis VCAs is as follows. Write the type (1) VCAs, $E(M^*,P^*)E(M^{*\prime},S^*)$: $S^*= S^*'M^*P^{*\prime}$, $P^*=S^{*\prime}M^{*\prime}P^*$, as $A(M^*,P^{*\prime})A(S^*,M^*)$: $A(S^*,P^{*\prime})$, where $S^* \in \{S, S'\}$, $M^* \in \{M, M'\}$, $P^* \in \{P,P'\}$. These are eight Barbara VCAs in the $S^*$, $M^*$, $P^{*\prime}$ variables. (The transition from one triplet of $S^*$, $M^*$, $P^{*\prime}$ values to another can be interpreted as a terms' relabeling transformation – for more details see below. But for now, as a concrete example, consider the VCA, $A(M',P)A(S',M')$: $A(S',P)$, and re-write it as a standard Barbara VCA – the one which may be written as $A(M,P)A(S,M)$: $A(S,P)$. The relabeling transforming one into the other is $ms(A(M',P)A(S',M'): A(S',P))=A(M,P)A(S,M): A(S,P)$, where $m(M)=M'$, $m(M')=M$, $s(S)=S'$, $s(S')=S$. Equivalently, denote $\mu:=M'$, $\sigma:=S'$; then the example becomes $A(\mu,P)A(\sigma,\mu): A(\sigma,P)$ which is a "standard Barbara" in the variables $\mu,\sigma,P$.) Still another way to see that any VCA of type (1) may be written as Barbara is to show that Barbara may be written as any other VCA of type (1): using the previous notation for $\mu$ and $\sigma$, Barbara's PCP and LC may be written as: $A(M,P)=A(\text{non } M',P)=A(\mu',P)$, $A(S,M)=A(\text{non } S',\text{non } M')=A(\sigma',\mu')$, $A(S,P)=A(\text{non } S',P)=A(\sigma',P)$. Thus Barbara can also be written as $A(\mu',P)A(\sigma',\mu'): A(\sigma',P)$. One can write the type (2) VCAs, $E(M^*,P^*)E(M^*,S^*)$: $M^*= S^{*\prime}M^*P^{*\prime}$, as $A(M^*,P^{*\prime}) A(M^*,S^{*\prime})$: $I(S^{*\prime},P^{*\prime})$ if $M^* \neq \emptyset$. These are eight Darapti VCAs in the $M^*$, $S^{*\prime}$, $P^{*\prime}$ variables. One can write the type (3a) VCAs, $E(M^*,P^*) I(M^*,S^*)$: $M^*S^*P^{*\prime}\neq\emptyset$, as $A(M^*,P^{*\prime}) I(M^*,S^*)$: $I(S^*,P^{*\prime})$. These are eight Darii VCAs in the M*, S*, P*' variables. (One could have kept the middle term in the LCs – the results would have been the same: the four VCA types can be written, respectively as a Barbara, Darapti, Darii and Disamis VCAs.) For the type (3b) Disamis VCAs, I(M*,P*) E(M*,S*): M*P*S*'≠Ø, switch the premises order and apply the convention that the firstly listed premise contains the P term, or, equivalently, perform the relettering S↔P, then switch the premises, to obtain exactly the eight Darii (3a) VCAs. Thus any VCA may be written as one of the Barbara, Darapti, Darii, or Datisi VCA. The detailed way to show that each of the four VCA types can be written, respectively as a Barbara, Darapti, Darii and Disamis VCAs is described below.

One can define a "relabeling group", G, acting on all 64 distinct PCPs of STM. Let p:= P↔P', s:=S↔S', m:=M↔M'. One can see that the compositions of s,p,m generate a commutative group G with eight distinct elements: 1,s,p,m, sp,sm,pm,spm. Obviously $1=s^2=p^2=m^2=(spm)^2=(ms)^2=(ps)^2=(pm)^2$. It is easy to see that each of the four types, (resp. subtypes), of VCAs, (resp. PCPs which do not entail any LCs), (1), (2), (3a), (3b), (resp. (4a), (4b), (5a), and (5b)), contains eight VCAs, (resp. PCPs), and that each type or subtype is left invariant under G. The invariance under G of each of the above types and subtypes is clear since each group, (of eight VCAs or PCPs), is described by a unique formula, containing S*, P*, M*, where S*∈{S,S'}, P*∈{P,P'}, M*∈{M,M'}, and that formula preserves its shape under G. (See formulas (α), (β), (γ), (δ) from Section 2, and formulas (4a), (4b), (5a), and (5b) toward the end of the Introduction.) Some examples of how the VCA types (1), (2), (3a), (3b) transform under G follow. For example, the m relabeling transforms Barbara into a Barbara', (All M' is P, All S is M', therefore S=SPM', etc.), and vice-versa; m also transforms Darapti into a Darapti' and vice-versa, m(Celarent/Cesare)=(Camestres/Camenes), etc. The p relabeling: p(Celarent/Cesare)= Barbara, etc. (The prime added to the name of a VCA shows that its LC refers to a subset of M' as non-empty or possibly non-empty, while the LS of the same VCA name without a prime refers to a corresponding subset of M: Barbara's LS is S=SPM, the LS of Barbara' is S=SPM', etc.) As another example, theVCA, E'E':M'=M'S'P', transforms under spm into the VCA Darapti: spm(E'E')= spm(E(M',P)E(M',S))=E(M,P')E(M,S')=A(M,P) A(M,S)=AA which are Darapti's premises, and spm(M'=M'S'P') becomes Darapti's LC, M=MSP. Equivalently, by contraposition, one may read E'E', i.e., All P is M, All S is M, as, All M' is P', All M' is S', which is now a Darapti in the variables M',P',S', with the LC M'=M'S'P', which is still the LC of E'E'. Note that there are no "free choices" in the term names: the term which appears in both premises is denoted M or M'; the term which appears in the firstly, (resp. secondly), listed premise is denoted P or P', (resp. S or S'). Since, via term relabelings, all variability inside each of the eight VCA sets (1), (2), (3a) and (3b) is accounted for, and all relabelings depend only on which set is denoted by either S or S', which one is denoted by either P or P', and which one is denoted by either M or M', it results that all eight VCAs from each of the four types of VCAs are logically equivalent. Same discussion applies to the eight PCPs from each of the subtypes (4a), (4b), (5a), and (5b) – but they continue to be uninteresting since they do not entail any LCs. The S↔P relettering, when combined with the convention to firstly list the P premise, amounts to changing the premises' order. CCS uses it for syllogism reduction to the 1$^{st}$ figure; it transforms the type (3a) VCAs into type (3b) VCAs and vice versa. To see precisely how G acts on the VCAs, one may arrange their 32 LC entailing PCPs into eight sets of four PCPs each, such that the variables in each set are (almost) the same. These eight sets, each containing four PCPs, will be denoted from 1 to 8, and one can let the relabeling group G act at once on each of the sets denoted by 1 to 8. Below, the first PCP in such a set of four PCPs always belongs to type (2) PCPs, the following two belong to the types (3a) and (3b) PCPs and the last one belongs to type (1) PCPs. If the precise LCs are kept, (pointing to just one subset of U), a PCP of type (1) – Barbara, entails, according to formula (α) of Section 2, two LCs and thus generates two VCAs, from which two ei VCAs can further be inferred. If M* is dropped from the LC, the above two VCA coincide, but the two ei VCAs will still be distinct. As formula (β) from Section 2 shows, any type (2) PCP generates one VCA with the LC M*=M*S'*P'* which will be discarded if one wants the M* eliminated, and only an ei VCA will be generated - via imposing the ei condition M*≠Ø. Thus, if the precise, one subset pointing to, LCs are kept, then the total count of the VCAs generated by the following eight groups of four PCPs each, comes to 40 (non-ei) VCAs and 24 ei VCAs. When the middle term M* is eliminated, then the total count of the VCAs generated by the following eight groups of four PCPs each, comes to 24 VCAs and 24 ei VCAs, since the number of the VCAs decreases by a total of 16, (a decrease of eight in the number of VCAs generated by each of the type (1) and type (2) PCPs). We'll say that each of such a four PCPs set, 1 to 8, is "bound to" the same subset of U: the four PCPs do not act at all on the subset of U on which the four PCPs are all "bound", but act on some of its "neighbors" in the Karnaugh map. Thus to each of the eight subsets of U one "attaches" a set of four PCPs "bound" to it. Listing on one column the four PCPs, and on the second column their LCs, these eight sets of four PCPs and and their entailed LCs are:

1. VCAs bound to the subset S'P'M:
EE=E(M,P)E(M,S)            M=S'P'M. If M≠Ø: I(S',P'), No names
IE=I(M,P)E(M,S)              S'PM≠Ø or O(P,S), No name
EI=E(M,P)I(M,S)              SP'M≠Ø or O(S,P), Ferio/Festino/Ferison/Fresison
EE'=E(M,P)E(M',S)           S=SP'M, P=S'PM', E(S,P),  Celarent/Cesare
                                       O(S,P) if S≠Ø, Celaront/Cesaro; O(P,S) if P≠Ø, No name

2. VCAs bound to the subset SP'M:
EA=E(M,P)E(M,S')             M=SP'M non-ei Felapton. If M≠Ø: O(S,P), the ei Felapton/Fesapo
IA=I(M,P)E(M,S')             SPM≠Ø or I(S,P), Disamis/Dimaris
EO=E(M,P)I(M,S')             S'P'M≠Ø or I(S',P'), No name
EA'=E(M,P)E(M',S')           P=SPM', S'=S'P'M, A(P,S)=A(S',P'),
                                I(S,P) if P≠Ø, Bramantip' (the prime refers to M' in
                                  P=SPM'); I(S',P') if S'≠Ø, No name

3. VCAs bound to the subset S'PM:
AE=E(M,P')E(M,S)             M=S'PM. If M≠Ø: O(P,S), No names
OE=I(M,P')E(M,S)             S'P'M≠Ø or I(S',P'), No name
AI=E(M,P')I(M,S)             SPM≠Ø or I(S,P), Darii/Datisi
AE'=E(M,P')E(M',S)           S=SPM, P'=S'P'M', A(S,P),  Barbara
                                 I(S,P) if S≠Ø, Barbari; I(S',P') if P'≠Ø, No name

4. VCAs bound to the subset SPM:
AA=E(M,P')E(M,S')            M=SPM – the non-ei Darapti. If M≠Ø: I(S,P), the ei Darapti
OA=I(M,P')E(M,S')             SP'M ≠Ø or  O(S,P), Bocardo
AO=E(M,P')I(M,S')             S'PM≠Ø or O(P,S), No name
AA'=E(M,P')E(M',S')           S'=S'PM, P'=SPM', E(S',P'), No name
                                 O(P,S) if S'≠Ø, No name; O(S,P) if P'≠Ø, No name

M' row VCAs:
5. VCAs bound to the subset S'P'M':
E'E'=E(M',P)E(M',S)          M'=S'P'M'. If M'≠Ø: I(S',P'), No names
I'E'=I(M',P)E(M',S)             S'PM'≠Ø or O(P,S), No name
E'I'=E(M',P)I(M',S)             SP'M'≠Ø or O(S,P), Baroco
E'E=E(M',P)E(M,S)              S=SP'M', P=S'PM, E(S,P),  Camestres/Camenes

O(S,P) if S≠∅, Camestros/Camenos; O(P,S) if P≠∅, No name

6. VCAs bound to the subset SP'M':
E'A'=E(M',P)E(M',S')         M'=SP'M' non-ei Felapton'. If M≠∅: O(S,P), the ei Felapton'/Fesapo'
I'A'=I(M',P)E(M',S')         SPM'≠∅ or I(S,P), Disamis'/Dimaris'
E'O'=E(M',P)I(M',S')         S'P'M≠∅ or I(S',P'), No name
E'A=E(M',P)E(M,S')           S'=S'P'M', P=SPM, E(S',P)=A(P,S), No name
                             I(S,P) if P≠∅, Bramantip, I(S',P') if S'≠∅, No name

7. VCAs bound to the subset S'PM':
A'E'=E(M',P')E(M',S)         M'=S'PM'. If M'≠∅: O(P,S), No names
O'E'=I(M',P')E(M',S)         S'P'M'≠∅ or I(S',P'), No name
A' I'=E(M',P')I(M',S)        SPM'≠∅ or I(S,P), Darii'/Datisi'
A'E=E(M',P')E(M,S)           S=SPM', P'=S'P'M, A(S,P)=A(P',S'), Barbara'
                             I(S,P) if S≠∅, Barbari'; I(S',P') if P'≠∅, No name

8. VCAs bound to the subset SPM':
A'A'=E(M',P')E(M',S')        M'=SPM' – the non-ei Darapti'. If M'≠∅: I(S,P), the ei Darapti'
O'A'=I(M',P')E(M',S')        SP'M' ≠∅ or O(S,P), Bocardo'
A'O'=E(M',P')I(M',S')        S'PM≠∅ or O(P,S), No name
A'A=E(M',P')E(M,S')          S'=S'PM, P'=SP'M, E(S',P'), No name
                             O(P,S) if S'≠∅, No name; O(S,P) if P'≠∅, No name

One can check that the group G transforms - at once - the above eight PCPs and LCs sets, denoted 1,2,...,8, as follows: s(i)=i+1, for i=1,3,5,7, and, since $s^2=1$, s(i)=i-1, for i=2,4,6,8; p(i)=i+2, for i=1,2,5,6, and, since $p^2=1$, p(i)=i-2, for i=3,4,7,8; m(i)=i+4, for m=1,2,3,4, and, since $m^2=1$, m(i)=i-4, for m=5,6,7,8. Except for E'A and A'A, the PCPs from groups 5 to 8 transform via an m relabeling into PCPs from groups 1 to 4 containing only positive terms. E'A contains only positive terms already, and A'A needs a p relabeling to become a positive terms only PCP. The PCPs EA' and AA' from groups 1 to 4, after an s relabeling, become PCPs containing only positive terms. Thus, since modulo a relabeling m:M↔M', almost all PCPs are formulable via only positive terms, one may say that CCS is justified to consider only PCPs with positive terms, and that neglecting PCPs whose LCs are A(P,S) or O(P,S) is also justified – a change in premises' order, which also implies an S↔P relettering, would transform those VCAs into VS.

One can also check that {G(1)}= {G(2)}= …={G(8)}= {1,2,3,...,8}={The set of all VCAs, with or without the middle term being eliminated}. (For example, {G(1)}={e(1)=1, s(1)=2, p(1)=3, sp(1)=4, m(1)=5, sm(1)=6, sp(1)=7, spm(1)=8}={The set of all VCAs, with or without the middle term being eliminated}.) This shows again that any VCA from any of the four VCA types can be recast as any other VCA of the same type. (A change of the order of premises, together with the convention of firstly listing the P-premise, completes the transformation of a type (3a) VCA into a type (3b) VCA and vice versa.)

The above represents a complete proof that any non-ei or ei VCA, (and VS), with the middle term eliminated from its LCs, is equivalent, modulo one of the G group relabelings, to either Barbara, the ei Barbari - provided one imposes the ei condition S≠∅, or another ei VCA with the same premises as Barbara and having I(S',P') as LC - provided one imposes the ei condition P'≠∅, or to a non-ei or ei Darapti, or to a Darii, or to a Disamis. If one does not take into account the ei VCAs, and one refers only to the non-ei VCAs whose LCs are the precise ones – out of which the middle term was not eliminated – one may say, in short, that any such VCA is equivalent to either Barbara – whose two LCs are S=SPM, and P'=S'P'M', or to a non-ei Darapti, or to a Darii, or to a Disamis.

## 8. About empty sets

The four types of VCAs may also be used to settle which VCAs are compatible with some of the sets S,P,M,S',P',M' being empty. In the modern square of opposition A(M,P), E(M,P) are not contraries anymore - unless one adds the condition M≠Ø. Instead when both A(M,P), E(M,P) are true it results that M=Ø. This empty set constraint (ESC) – which empties all four subsets of M- is compatible with the universal premises of the VCAs of type 1 and 2 - but not with the ei on M. Nor is the M=Ø ESC compatible with the type (3a) and (3b) premises, I(S*, M) and I(P*, M). In fact, since the VCAs of the same type are all equivalent, it results that a complete discussion of the compatibility of various ESCs and VCAs may be reduced to examining just three or four representative cases. Moreover, instead of firstly imposing an ESC, and then finding out the PCPs compatible with it, one can do it the other way around, by listing, for each VCA type, the ESCs with which that VCA type is compatible or incompatible.

Darii's PCP, A(M,P)I(S,M), means MP'=Ø, SM≠Ø, and the LC is SM= SMP+ SMP'= SMP≠Ø. From the LC SMP≠Ø, one may, with some loss of information, eliminate M, and re-express the LC as I(S,P)="Some S is P". Thus the PCP is incompatible with the S=Ø, M=Ø, and P=Ø ESCs, but is compatible with the S'=M'=P'=Ø ESCs, (which imply S=M=P=U; thus in this latter, extreme, case Darii's PCP and LC just assert that U is non-empty).

Darapti's PCP, A(M,P)A(M,S), means MP'=Ø, MS'=Ø, and the LC is M=MP+ MP'=MP= MPS+MPS'=MPS, which may be written as A(M,SP). This time around one may eliminate M only via the ei hypothesis M≠Ø, then re-express the LC as I(S,P). Thus the ei hypothesis is incompatible with the M=Ø, S=Ø and P=Ø ESCs, but is compatible with the S'=M'= P'=Ø ESCs, (which imply S=M=P=U; therefore, in this latter, extreme case, Darapti's PCP plus ei LC affirm that U is non-empty). Note that Darapti's PCP without the added ei condition is compatible even with U=Ø, in which case the PCP is just "chatter about empty sets".

Barbara's PCP, A(M,P)A(S,M), means MP'=Ø, SM'=Ø, and the LCs are S=SM+ SM'=SM=SMP+ SMP'=SMP, and P'=P'M+ P'M'=P'M'= P'M'S+P'M'S'=P'M'S'. The first LC may be written as A(S,MP), or, with some loss of information, one may eliminate M, and write A(S,P) =E(S,P'), which now refers to an entire column of U instead of just one of the eight subsets of U. (A "precise" LC always pinpoints to just one of the eight subsets of U.) The second LC may be written as A(P', S'M'), or, with some loss of information, one may eliminate M', and write A(P',S') =E(S,P') – the same as the first LC. Since Barbara's PCP contains only universal premises, the PCP is compatible even with U= Ø in which case all the deductions and the LCs – either "precise" or "classically expressed", are just "chatter about empty sets". One may then add an ei hypothesis, S≠Ø, to the 1st LC, and a different ei hypothesis, P'≠Ø, to the 2nd LC, to obtain, after the M, resp. M', elimination the new ei LCs: I(S,P), (Barbari), and resp., (the un-named), I(S',P'). The S≠Ø ei hypothesis means, since S=SPM, that also P≠Ø and M≠Ø, while the compatible ESCs are S'= Ø, or/and, P'= Ø, or/and, M'=Ø. The S'=P'=M'= Ø constraint amounts to Barbari affirming U≠Ø. The P'≠Ø ei hypothesis means that also S'≠Ø and M'≠Ø. If both ei hypotheses are true then all the sets M,M',S,S',P,P' are non-empty, and there are no ESCs compatible with both ei hypotheses.

In conclusion any universal premise is compatible with any ESC. But any ei hypothesis or any LC of a VCA of type (3a) or (3b), (containing one universal and one particular premise - both acting on either M or M'), specifies three sets that are non-empty, and thus pinpoints to three ESCs with which the ei hypothesis or the "type (3a) or (3b) LC" is incompatible. The above considerations are based on a sort of "temporal commutativity": instead of firstly applying the ESC to obtain a particular universe of discourse, and then searching for the LC in that universe, one writes down the LC in the usual 8-subset universe of discourse U, and one applies the ESC only afterwards, to see if it is compatible with the PCP and its LC.

## 9. An expanded definition of distribution, and the RofVCA regarding distribution conservation

"A term is said to be distributed when reference is made to all the individuals denoted by it; it is said to be undistributed when they are only referred to partially, i.e., information is given with regard to a portion of the class denoted by the term, but we are left in ignorance with regard to the remainder of the class." (Keynes, p.68). Firstly, one may expand the definition of distribution, by agreeing that whatever distribution the two terms appearing in a categorical statement may have, then their complementary sets/terms in U, are automatically assigned an opposite distribution. Thus since in I(M,P), the terms M,P are undistributed, the terms M',P' are distributed in the same I(M,P) statement. This is in agreement with the obversion and contraposition rules and the standard definition of distribution: from I(M,P)=O(M,P')= O(P,M') one realizes that M' and P' are distributed since M,P were not. This expanded definition of the distribution allows us to notice that, e.g., in a Darapti VCA, (and VS), E(M,P)E(M,S): A(M,S'∩P'); I(S',P') if M≠Ø, the distribution of M is the same in the premises as it is in the precise LC, and the distributions of S and P are the same in the premises as they are in the I(S',P') LC, after M was eliminated: S',P' are undistributed in the LC, as they were in Darapti's premises, since S,P were distributed there. Thus the distribution of the middle term is conserved in the precise LC, and the distributions of the end terms are conserved in the LC from which the middle term was eliminated. The above expresses a generalization to all the VCAs, (even to the ei VCAs of type (2)), of the usual RofVS - which refer to only non-ei VS: "if a term is distributed in the LC, it must be distributed in its corresponding premise". As the formulas (α),(β),(γ),(δ) from Section 2 prove it by inspection, the following is true:
RofVCA #1. The distribution of the end terms is conserved in all non-ei and ei VCAs, except in type (1) ei VCAs, where ei on S*, (resp. P*) changes S*, (resp. P*), from distributed in the PCP to undistributed in the ei LC, while the distribution of the other end term, P*, (resp. S*), remains the same as it was in the PCP.
Note also that the arguments of a statement and of its contradictory one, have opposite distributions – in E(M,P) both M and P are distributed, while in the contradictory statement, I(M,P), both M and P are undistributed, (while M',P' are distributed). Similarly, in A(M,P), M is distributed and P is not, while in the contradictory statement, O(M,P), the terms' distributions are reversed.

## 10. A proof of RofVCA #4 - "if one premise is negative, the conclusion must be negative, if the premises are both affirmative or both negative, the conclusion is affirmative"

Firstly one needs to define, for both universal and particular premises, when they are considered as affirmative and when they are considered as negative statements. Then, one has to show that RofVCA #4 is satisfied for each choice of S*,M*,P*, in each of the formulas (α), (β), (γ), (δ). The universal negative premises are E(M,P), E(M,S), E(M',P'), E(M',S'), and the only particular negative premises are O(P,M), O(S,M), O(M,P), and O(M,S). Denoting h∈ {S,P}, h'∈ {S',P'}, the universal negative premises are E(M,h),

E(M',h'), the particular negative premises are I(M',h), I(M,h'), the universal affirmative premises are E(M',h)=A(h,M), E(M,h')=A(M,h), the particular affirmative premises are I(M,h), I(M',h'). One can see that the switch E↔I while the arguments are left unchanged transforms universal negative premises into particular affirmative premises and vice versa, and transforms universal affirmative premises into particular negative premises and vice versa. The switch M↔M', (resp. h↔h'), transforms affirmative premises into negative premises and vice versa. Based on these observations, and on the formulas α, β, γ, δ, one can verify RofVCA #4: Two affirmative premises, (as in Barbara or Darii), entail, after middle term elimination, an affirmative LC; two negative premises, (as in Barbara' or Darii', i.e., same premises as Barbara and Darii except for the relabeling M↔M'; see Section 7), entail, after middle term elimination, an affirmative LC; one affirmative plus one negative premises entail, after middle term elimination, a negative LC - as in Celarent/Cesare, (aka Camestres'/ Camenes'), and Camestres/Camenes, (aka Celarent'/Cesare'). For example: E(M,P) E(S,M')=E(M,P) A(S,M) – which is the Celarent/Cesare's PCP, entails two LCs, which, after the middle term elimination, become the same, negative, and universal statement E(S,P). But two negative premises, E(M,P)E(M,S), entail a positive ei LC: M=MP'S', M≠Ø, so I(S',P'). (One sees, once more, that the positive and negative terms appear "on equal footing". One can prove the RofVCA #4 by defining the negativity or signature, s, of a statement symbol, s(A)=s(I)=0, s(E)=s(O)=1, the signature of a term, s(M)=s(P)=s(S)=0, s(M')=s(P')=s(S')=1, and, the signature of a whole statement as the sum of the signatures modulo 2 of the statement's symbol and all of its terms. Then a statement is affirmative, (or positive), if its signature is zero, and is negative if its signature is 1. Thus, e.g., s(A(M,P)) =s(E(M,P'))= s(A(P',M'))=0, s(A(M',P))=s(E(M',P'))=1; therefore the first three statements are affirmative, and the last two are negative. Now one can prove by cases the above RofVCA #4 for each of the four VCA types. For type (1) Barbara, with premises E(M*,P*)E(M*',S*), and with the LCs S*=S*M*P*', P*=P*M*S*' or A(S*,P*')= A(P*,S*')=E(S*,P*), the premises' signatures are (1+s(M*)+s(P*)) mod 2, (2+s(M*)+s(S*)) mod 2, and the LC, after middle term elimination, has the signature (1+s(P*)+s(S*)) mod 2. It results that when s(S*)=s(P*) the LC is a negative statement, while the premises have different signatures, i.e., one premise is affirmative and one negative; thus RofVCA #4 holds. When s(S*) and s(P*) have different signatures then the LC is an affirmative statement, while the premises have the same signature, i.e., the premises are either both affirmative or both negative – and RofVCA #4 holds again.

Instead of continuing to verify RofVCA #4 in the same way for the formulas (β), (γ), (δ), one can slightly simplify the verification procedure by observing that the formulas (α), (β), (γ), (δ) - from which the middle term was eliminated - were already written to formally satisfy all the RofVCA, including RofVCA #4: For types (1) and (2) – the two premises are negative, E(M*,P*)E(M*',S*), (resp. E(M*,P*)E(M*,S*)), and, according to the RofVCA #4, the LCs are affirmative: A(S*,P*'), (resp. I(S*',P*')). (In these affirmative LCs the distributions of S* and P*are exactly the same as they were in the PCPs – and this proves RofVCA #1 for PCPs of types (1) and (2).) Premises of types (3a) and (3b), (one negative and one affirmative premises), E(M*,P*) I(M*,S*), (resp. I(M*,P*) E(M*,S*)), entail the negative LCs, O(S*,P*), (resp. O(P*,S*)), in which the distributions of S* and P* are again the same as they were in the PCPs. To prove RofVCA #4 in a simpler way, it remains to show that if the terms signatures in the LCs change the LC statement's signature, those terms signatures change the premises signatures in such a way that the RofVCA #4 is still satisfied. Therefore, this slightly simpler proof of RofVCA #4 is based on dropping the statement symbol signatures, and considering only the statements' arguments signatures, via introducing the "argumental signature" of each premise, and of the LC. For type (1) Barbara, the argumental signature of each negative premise, E(M*,P*), (resp. E(M*',S*)), is [s(M*)+s(P*)] mod 2, (resp. [s(M*)+1+s(S*)] mod 2), and the argumental signature of the LC, A(S*,P*'), is [1+s(P*)+s(S*)] mod 2. It follows that when s(P*)=s(S*) the LC is in fact negative, but then the premises are one affirmative and one negative – as they should be according the RofVCA #4. When s(P*)≠s(S*), the LC remains affirmative, and the two premises have the same argumental signatures – thus the two premises will both remain negative, or both become affirmative. For the other formulas, β, γ, δ, the proof of the RofVCA #4 proceeds similarly: For type (2)

Darapti, E(M*,P*) E(M*,S*) with the LC I(P*',S*') if M≠Ø, the premises' argumental signatures are (s(M*)+s(P*)) mod 2, (s(M*)+s(S*)) mod 2, and the LC has the argumental signature (2+s(P*)+s(S*)) mod 2. It results that when s(S*)=s(P*), the LC remains an affirmative statement, while the premises have the same signature, i.e., both premises are affirmative or both are negative. When s(S*)≠s(P*) the LC becomes a negative statement, while the premises will have different signatures, i.e., one premise is affirmative and the other is negative. For type (3a) Darii, E(M*,P*) I(M*,S*) with the LC O(S*,P*), the premises' argumental signatures are (s(M*)+s(P*)) mod 2, (s(M*)+s(S*)) mod 2, and the LC has the argumental signature (s(P*)+s(S*)) mod 2. It results that when s(S*)≠s(P*) the LC becomes an affirmative statement, while the premises will have different argumental signatures, i.e., both premises are either affirmative or negative. When s(S*)=s(P*) then the LC remains a negative statement, while the premises have the same argumental signatures, i.e., the premises remain one affirmative and one negative. The proof that the type (3b) VCAs also satisfy the RofVCA #4 is identical to the one for type (3a) VCAs. One has thus checked, too, that the relabelings from the G group are in agreement – or "covariant" - with the RofVCA #4: any relabeling from G, when applied to any instance of one of the formulas (α), (β), (γ), (δ), transforms it into another instance of the same formula, such that both instances satisfy the RofVCA #4.

## 11.  An analysis of the RofVS "theory"

First thing to note, is that the RofVS, like the CCS, suppose that the PCPs are formulable using only positive terms, and that, if a PCP entails an LC, then the LC is formulable using only positive terms, too. One lists the RofVS, (Copi 2009, Hurley 2008, Stebbing 1961, Keynes 1887), in the following order: #1 - "the middle term has to be distributed in at least one premise", #2 - "two negative premises are not allowed", (these two syllogistic rules refer to the PCPs only – the following four rules refer to the entire syllogism, PCP and LC, and will be used to predict the LCs of any PCP from the RofVS' DofA), #3 - "any term distributed in the LC must be distributed in the PCP", (i.e., the distribution of the end terms, P and S can not "increase" from undistributed in the premises to distributed in the LC, but can, conceivable, decrease from distributed in the premises to undistributed in the LC), #4 - "if either premise is negative, the LC must be negative", #5 - "from two universal premises, no particular LC may be drawn", #6 - "if one premise is particular, then the LC is particular". The RofVS were developed into a more rounded "theory": firstly, one specifies a set of PCPs which is not part of the RofVS' DofA: according to RofVS #1, "PCPs in which the middle term is not distributed in any of the two premises, do not entail any LCs", and, according to RofVS #2, "PCPs made of two negative premises do not entail any LCs". (Copi 2009, Hurley 2008, show that such PCPs, indeed, can not entail LCs formulable via only positive terms, but Carroll, 1977, p.240, and the formulas (α), (β), (γ), (δ), show that such PCPs can entail LCs of the format I(S',P'), thus generating VCAs.) After the set of PCPs to which the RofVS do not apply is postulated by RofVS #1 and #2, one can deduce that there are other PCP sets for which the RofVS will predict LCs which contradict the RofVS themselves: therefore, those PCP sets will not be part of the RofVS' DofA, either. What is left from the set of PCPs formulable via only positive terms, after the above two rounds of PCP removals, should be mostly those PCPs on which the CCS' and the RofVS' results coincide: the PCPs which generate VS, and the PCPs whose entailed LCs are A(P,S) and O(P,S) – the latter PCPs also satisfy all the RofVS, as the VS do. As Stebbing 1961 pp. 86-92, shows, no PCPs made of two particular premises, may produce LCs compatible with the RofVS. Indeed, PCPs made of two particular and affirmative premises – in which, therefore, no term is distributed, (resp. two particular and negative premises – the "not allowed" PCPs), are already excluded from the DofA, by RofVS #1, (resp. RofVS #2), and the PCPs made of one affirmative and one negative particular premises, should have an LC which is particular and negative, (according to RofVS #4 and #6). But this means that an end term will be distributed in the LC, without being distributed in the PCP – since, according to the distribution's definition, only the middle term will be distributed in the negative particular premise: contradiction. Thus, by postulating two new classes of PCPs which do not entail LCs, one was able to prove that any PCP made of two particular premises does not entail an LC – therefore all PCPs of types (4a) and (4b) will be excluded from the DofA of the RofVS. One may now prove that PCPs of types (5a) and (5b) do not entail any LCs, either. Namely, one can check that out of the four PCPs of type (5a) and four PCPs of type (5b) which are formulable only via positive terms, two of them contain only negative premises, in four of the PCPs the middle term is not distributed at all, and in another two PCPs one of the premises is particular, one is negative, both end terms are undistributed, and therefore, the LC should be negative and particular, thus one end term would be distributed in the LC, without being distributed in the premises - thus contradicting RofVS #3. [The complete details are as follows. The formula (5a), $E(M^{*'},P^*)I(M^*,S^*)$, leads, for $M^*=M$, $P^*=P$ and $S^* \in \{S, S'\}$ to two PCPs, (containing only positive terms), where M is nowhere distributed – thus, in accordance with RofVS #1 these PCPs do not entail an LC. Similarly, the formula (5b), $I(M^*,P^*) E(M^{*'},S^*)$, leads, for $M^*=M$, $S^*=S$ and $P^* \in \{P, P'\}$ to two PCPs, (containing only positive terms), where M is nowhere distributed. According to RofVS #1, all the above four PCPs will not entail any LC. In $A(M,P)O(S,M)$ – which is formula (5a) for $M^*=M'$, $S^*=S$ and $P^*=P'$, and in $O(P,M)A(M,S)$ – which is formula (5b) for $M^*=M'$, $S^*=S'$ and $P^*=P$, the middle term is distributed in both premises, but S and P are nowhere distributed; since the LC should be particular and negative, (due to one premise being negative and particular), the LC would distribute either S or P, thus contradicting RofVS #3. Finally two PCPs, $E(M,P)O(S,M)$ – formula (5a) for $M^*=M'$, $S^*=S$ and $P^*=P$, and $O(P,M) E(M,S)$ – formula (5b) for $M^*=M'$, $S^*=S$ and $P^*=P$, are made only of negative premises, and therefore they are not contained in the RofVS' DofA, in accordance with RofVS #2.]

After the RofVS' DofA was found as being made of the initial 36 PCPs formulable via only positive terms, out of which one removes the three PCPs where both premises are negative and at least one premise is universal, the one PCP containing two universal premises in which M is nowhere distributed, the nine PCPs containing only particular premises, and the eight PCPs of types (5a) and (5b), one can check that the RofVS #3 to #6 will predict, for most of the 36-(3+1+9+8)=15 PCPs left as candidates for RofVS' DofA, one of these LCs: A(P,S), O(P,S), A(S,P), E(S,P), I(S,P), O(S,P). One can check that the VCA generated in accordance to the RofVS, are either the eight standard Boolean VS, (described by CCS as 15 VS when the syllogistic figures are taken into account), or four VCAs having A(P,S) or O(P,S) as LCs, and which become VS via a metathesis and the relettering S↔P. (One may call them Bramanta, (which by ei on P gives the ei Bramantip), Fireo, Boraco and Bacordo. They satisfy all the RofVS, with the LCs, as predicted by the RofVS, being either A(P,S) – for Bramanta, or O(P,S) – for the other three.) These four latter VCAs can be discarded, exactly as A(P,M)A(M,S) - "Bramanta" - is discarded in CCS, on "relettering grounds". The three PCPs left, (out of the 15 PCPs constituting the RofVS' DofA candidates), can not be handled by the RofVS – since the RofVS do not refer at all to the ei VS, but the RofVCA – which predict LCs for both VCAs and ei VCAs, predict ei LCs for these three PCPs – namely they will generate the ei Darapti, Felapton/Fesapo, and "Falepton" - the latter having O(P,S) as LC. For example, the presentations of the RofVS in Copi and Hurley start by showing that all the (eight standard Boolean) VS identified by CCS, necessarily satisfy all the RofVS - see Hurley (pp. 285-289), Copi (pp.230-235). Proceeding in the reverse direction, one shows that, the only PCPs whose possible entailed LCs are A(S,P), E(S,P), I(S,P), O(S,P), and such that the PCPs and their entailed LCs satisfy all the RofVS, are the eight standard Boolean VS, i.e., the 15 standard (non-ei) VS described by CCS - see Hurley (pp.290-291), Copi (pp.245-248). Nevertheless this does not show that the only PCPs which satisfy all the RofVS are the eight standard Boolean VS, (Bramanta, Fireo, Boraco and Bacordo also satisfy all the RofVS!) – it shows that the only PCPs which generate VS and satisfy all the RofVS, are the eight standard Boolean VS described by CCS. (Reminder: a VS is generated by a PCP formulable via only positive terms, which validly entails an LC of one of the A(S,P), E(S,P), I(S,P), O(S,P) formats. And these are the only LC formats considered by the proofs in Hurley (pp.290-291) and Copi (pp.245-248).)

## 12. An analysis of the RofVCA "theory"

The RofVCA "theory" can be presented in a similar way to the RofVS "theory". One defines the Domain of Applicability (DofA) for the RofVCA #1 to #4, as the maximal set on which the STM – which produced the formulas α, β, γ, δ, can be practically replaced by the RofVCA #1 to #4, because the latter predict the same LCs as the ones the STM proved are true. Firstly, one postulates, (and this is a true fact -  which can not be said about the postulates expressed by the RofVS #1 and #2), that PCPs of type (4a) are not part of DofA for the RofVCA. Then one can prove that PCPs of types (4b), (5a) and (5b) are not part of DofA for the RofVCA, either, since when applied to such PCPs, the RofVCA will predict LCs which contradict the RofVCA themselves. The PCPs left, are those of types (1), (2), (3a) and (3b), on which – as proved in previous Sections - the results of RofVCA and STM coincide. For example, applying the RofVCA to (5a) PCPs, $E(M^{*'},P^*)I(M^*,S^*)$, one obtains – since one premise is negative and one is affirmative and particular – that the LC is $O(S^*,P^*)$. Therefore if $s(S^*)=s(P^*)$, (resp. $s(S^*){\neq}s(P^*)$), the LC remains a negative statement, (resp. becomes an affirmative statement), and the premises will have opposite, (resp. identical), argumental signatures – this way the premises become either both negative or both affirmative, (resp. the premises remain one negative and one affirmative), contrary to what the RofVCA #4 said – while explicitly being used to find the false LC. (As one already knows from the Introduction, the definition of a syllogism implies that PCPs of types (4a), (4b), (5a) and (5b) do not entail any LCs, anyhow.) If one also extends the RofVCA #3 to say that if at least one premise is particular then the LC is particular, then these "extended" RofVCA applied to PCPs of type (4b), would predict LCs which are contradicted by the RofVCA themselves. The proof is similar to the one given above for PCPs of type (5a). The proof that for PCPs of type (5b),  the RofVCA will predict LCs which contradict the RofVCA themselves – is also similar to the proof given above that PCPs of type (5a) are not part of DofA for the RofVCA.

Consider the PCPs E'E', EE, OE, EO - whose LCs, (or ei LCs), as one already knows, are all the same: $I(S',P')$. One can see that the RofVCA predict the correct LC for these PCPs: In the case of $E'E'=E(M',P)E(M',S)=A(P,M)A(S,M)$, one has two affirmative premises, (as one knows, of type (2) Darapti). According to the RofVCA, (which do not depend on the PCP type, provided that the PCP is contained in the RofVCA' DofA), the LC should be affirmative and should conserve the fact that both S and P are distributed in the premises. But the only statement out of the eight possible LCs remaining after the middle term is eliminated, $E(S^*,P^*)$, $I(S^*,P^*)$, where $S^*\in\{S,S'\}$, $P^*\in\{P,P'\}$, which is affirmative and in which both S and P are distributed, is $I(S',P')$. If, instead, as it was done in Section 10, one counts only the statements signatures, and one applies the of RofVCA #4 "directly" to $E'E'=E(M',P)E(M',S)$, counting that the LC of two negative premises has to be an affirmative LC which conserves S and P being distributed in the premises, one again obtains that the unique LC which satisfies all the RofVCA #3 to #6, is $I(S',P')$. Indeed, trying instead $A(P,S')$ as LC of the above premises does not work, since $A(P,S') =E(P,S)$ which contradicts RofVCA #4. Thus, for the Darapti's PCP, the RofVCA predict even that the only VCA from which M is eliminated is an ei VCA. One sees once more, that on the set of PCPs of types (1), (2), (3a), and (3b), the RofVCA may effectively replace logic theory by "prediction rules". In cases, such as OE, EO, the RofVCA will similarly choose the same LC, $I(S',P')$ - based on the fact that two negative premises imply an affirmative LC, that both S and P are distributed in the premises, and that if one premise is particular,  the LC will be particular. If both premises are universal - as in EE -  the RofVCA, together, will still predict, as above, that the LC will be particular, (and thus, the ei condition – on M - is necessary).

## 13. How many sound VS or VCAs may one hope to construct out of three given terms, without imposing restrictions on the structure of the universal set U

When three specific terms are given, with one of them already designated as the middle term, one may consider all the 36 or 64 PCPs which can be constructed starting with these three specific terms, (out of which one is the designated middle term), and one can try to see what sound VS or VCAs one may construct out of the three terms.

As one shows below, given three terms, with one of them already designated as the middle term, then, at most one sound VCA of type (1) or (2) may be built out of the three terms without restricting U to particular cases. (This is also almost evident on the Figure 1 – the Karnaugh map - from Section 3.) Since that VCA can be formulated either as a Barbara or a Darapti if the terms are appropriately labeled, one may say that given three terms, there exists at most one sound VCA of types (1) and (2) – either a Barbara or a Darapti – which can be constructed out of the given three terms, (again, if one of them was already designated as the middle term). If the three given terms generate, (modulo a relabeling from G), a sound Barbara, then a maximum of two other type (3a) and two other type (3b) sound VCAs may perhaps be constructed with the same given three terms without restricting U to particular cases: these new sound VCAs have their universal premises "stolen" from Barbara and their possible particular premises place set elements on the four subsets adjacent to the four subsets emptied by Barbara's two universal premises. One of these other possible four VCAs is a Darii/Datisi, and the other three have no names since they assert that subsets – other than the three subsets preferred by CCS – are non-empty. If the three given terms generate, (modulo a relabeling from G), a sound Darapti, then only two other sound VCAs, one of type (3a) and one of type (3b) may be constructed with the same given three terms without restricting U to particular cases: these two new VCAs, a Darii/Datisi and a Disamis/Dimaris will have their universal premises "stolen" from Darapti, and the same LC as the ei Darapti (after the middle term is eliminated): $SPM \neq \emptyset$.

When the middle term is pre-determined, then two distinct PCPs of either type (1) Barbara, or type (2) Darapti, or one PCP of type (1) and one PCP of type (2), will necessarily contain either two distinct two P-premises, or two distinct S-premises, which will impose a particular structure on the universal set U. For example, if A(M,P) and E(M,P), are both true, as P-premises in two different PCPs, would imply M being empty, $M = \emptyset$. The relationships implied by the other five possible combinations of two universal P-premises being simultaneously true: E&E' imply $P = \emptyset$, A'&E' imply $M' = \emptyset$, A&A' imply $P' = \emptyset$, A&E' imply P=M, A'&E imply P=M', and similar relationships hold for the "top face" of the S-cube, (if one places the four universal premises on the top faces of the P and S cubes, and the particular premises on the bottom faces of the two cubes). This shows, e.g., that Barbara, AE', and Camestres, E'E, can both be sound, but uninteresting, since in that universal set, P=M and $S = \emptyset$.

Note that choosing another of the three terms as a middle term, leads to arguments and conclusions similar to the ones above: no two sound and distinct VCA of either type (1), or type (2) may be constructed with the same middle term, unless the universal set has a particular structure; no sound pair of a VCA of type (1) and one VCA of type (2) may be constructed with the same middle term, unless the universal set has a particular structure. Thus, the next task is to see if by using a term once as M, and a second time, say, as P, while the term firstly used as P, is afterward used as M, one can produce two distinct PCPs of type (1) without imposing, when all four premises are true, a particular structure on U. By adjoining the "standard" Barbara's PCP, A(M,P)A(S,M)=E(M,P')E(M',S), to each of the eight PCPs of type (1), having P as the middle term and S and M as end terms, E(P*,M*)E(P'*,S*), one can show that, given three terms, at most one sound VCA of type (1) may be constructed with them, without imposing a particular structure on the universal set U, i.e., in short, one may say, that at most one of the terms, (out of three given terms), may be used as the middle term in a type (1) VCA. For example, if all the following four premises, A(M,P)A(S,M) and A(P,M)A(S,P), (where the second Barbara's PCP is obtained by switching the roles which M and P played in the first Barbara's PCP), are true, i.e., MP'=SM'=PM'=SP'=0, then, M=MP'+MP=MP= MP+PM'=P. For a complete proof, one may compare, two at a time, the eight PCPs of type (1) having P as the middle term, with the "standard" Barbara's PCP. From E(M,P')E(M',S) and E(P,M) E(P',S*), it results M=0; from E(M,P')E(M',S) and E(P,M') E(P',S*), it results M=MP=P; from E(M,P')E(M',S) and E(P',M) E(P,S*), it results, if S*=S, that S=PS+P'S=P'S=P'SM+P'SM'=0+0=0, and, if S*=S', that M=MP=MPS=S; from E(M,P')E(M',S) and E(P',M') E(P,S*), it results P'=0.

One may easily see that the "standard" Barbara's PCP, E(M,P')E(M',S), and this type (2) PCP having S as the middle term, E(S,P')E(S,M'), can be simultaneously sound without imposing a particular structure on U. Even three type (2) PCPs can be simultaneously sound without imposing a particular structure on U, under the condition that their middle terms are all different. For example, E(M,P')E(M,S), E(S,P')E(S,M), E(P',M)E(P',S), use M, S, and resp. P' as middle terms, and out of their six premises only three are distinct. These three type (2) PCPs, empty a total of only four subsets of U – the same number as a single PCP of type (1) empties.

### 14. Why, despite their logical weak traits, the CCS and the RofVS were considered, for centuries, "good enough theories"?

In a few words, the answer seems to be – the relabelings from G, and the role reassignment S↔P. All inadvertencies and oversights that the CCS and the RofVs might have committed, were without consequences, because no matter if the premises, (and LCs), contain only positive terms or not, no matter if one accepts other LC formats except the A(S,P), E(S,P), I(S,P), O(S,P) – the ones prescribed by the VS definition, no matter if one discards VCA whose premises are both negative, or whose LCs are A(P,S) or O(P,S), no matter if one tries to find the VCAs which satisfy all the RofVS, while at the same time one supposes, from the very beginning, that their LCs have to be just those defining a VS, in the end, it did not matter, because there are only four distinct VCAs – Barbara, Darii, Disamis and the non-ei Darapti.

By comparison to STM and the RofVCA, the CCS and the RofVS seem contrived and not well put together. George Boole did not like the contrivance of the moods and syllogistic figures, nor the foreordained requirement that the term listed in the first premise has also to fulfil the role of being the predicate of the conclusion. (See Burris.) (I confess that the same uneasy feelings towards figures and P being the predicate of the LC, tempted me to think about terms as sets, and a possible STM, even before I finished my first reading through the syllogistic moods and the RofVS.)

In short, a partial list of the "sins" of CCS and RofVS seems to be:
- disregard premises which explicitly contain negative terms
- disregard LCs which explicitly contain negative terms
- not inventing a better description or characterization of the facts that A(M,P)≠A(P,M) and O(M,P)≠O(P,M) than the syllogistic figures
- decreeing not only that the term denoted by P has to appear in the first premise - good convention, but also that P has to be the predicate of the LC (this eliminates from consideration VCAs having A(P,S) and O(P,S) as LCs, while such VCAs satisfy all the RofVS, and while the syllogistic figures were invented exactly to describe such statements. Moreover, how one is supposed to determine the predicate of the LC – so that one could firstly lists that "major premise" - even before having any ideea if the PCP under examination does entail an LC?)
- not explicitly recognizing PCPs of types (5a) and (5b) as not entailing any LCs
- falsely claiming that a PCP which contains two negative premises never ever entails an LC

- falsely claiming that any PCP in which the middle term is nowhere distributed does not entail any LC
- the RofVS do not predict the ei LC of type (2), Darapti, PCPs
- obsessively eliminating the middle term from the LC, thus weakening it
- the STM "LC paradigm" - which states that an LC means a PCP pinpointing to a unique subset of U (the universe of discourse), never made it into CCS (and, up to now, seemingly not even into STM).

As mentioned above, all the above weaknesses never lead to a contradiction or a manifestly wrong prediction coming out of the CCS or the RofVS, mainly because the only four non-equivalent VS and VCAs are Barbara, non-ei Darapti, Darii, and Disamis.

## 15. Conclusions

The set model shows very intuitively what an LC is, and what a lack of LC means. STM, clearly defines the set of PCPs which generate VCAs. STM proves, via the formulas $(\alpha)$, $(\beta)$, $(\gamma)$, $(\delta)$, a generalized form of RofVS #3: that in a VCA, any of the end terms, S,S',P,P', is distributed in a premise if and only if it is distributed in the LC, unless ei is imposed on that term, in which case the term was distributed in its premise and is undistributed in the LC. STM also proves a generalized form of RofVS #4: in any VCA, either two affirmative, or two negative premises entail an affirmative LC; if the PCP contains only one negative premise, then the LC is also negative.

For classification purposes, the premise which contains a term denoted by P, is listed firstly in both STM and the CCS. The latter recognizes two equivalence classes of VCAs: the VS having A(S,P) or O(S,P) statements as LCs, and the VCAs with A(P,S) or O(P,S) LCs: the VCAs from one class are transformed into the other class VCAs via a relettering S↔P and/or a switch of the premises' order. This equivalence is taken as a sign that the VCAs from the second class may be discarded since they can be recasted as VS. Thus the VCAs whose LCs are A(P,S) or O(P,S) statements are weeded out by the very definition of a VS, not by RofVS, which the discarded VCAs still fully satisfy. Note though that the same S↔P relettering and a switch of the premises' order, transforms Camestres/Camenes into Celarent/Cesare, but nevertheless neither Camestres/Camenes nor Celarent/Cesare are discarded since both have as LC the E(S,P)= E(P,S) statement in which either term may be used as the predicate. This may suggest that the first two RofVS are rather ad hoc rules meant to weed out VCAs which contain negative terms in their PCPs or the LCs. The STM recognizes two similar classes of VCAs – namely the type (3a) and type (3b) VCAs. Via a metathesis and an S↔P relettering any of the two types may be rewritten as VCAs of the other type. Moreover any two VCAs of the same type, (either (1), (2), (3a) or (3b)), are equivalent via a relabeling from the 8-element group G, generated by the p:= P↔P', s:=S↔S', and m:=M↔M' relabelings and their compositions, ps, pm, sm, psm. Generalizing from VS to VCAs and from 36 PCPs to 64 PCPs, one realizes that modulo term relabelings there are only four distinct VCAs (and VS) – one per each of the four PCP types already discussed by George Boole and Lewis Carroll. One may try to use the structures of the sorites discussed in Section 4, to solve, e.g., the many sorite examples given by Carroll 1977, pp. 385-421. For newer takes on categorical syllogisms one may read, e.g., Frank Thomas Sautter, (2019), or Stephen Read, (2017). But the Boole and Carroll "four (or three) VCA types approach" so clearly displays the structure of any VCA, that, I think, after two millenia, one should take away the center stage from the VS moods and figures, and the first two RofVS.